\documentclass{amsart}
\usepackage{amssymb}
\usepackage{amscd}

\title[$p$-adic multiple zeta values II]
{$p$-adic multiple zeta values II
--- tannakian interpretations}
\author{Hidekazu Furusho}

\address{Graduate School of Mathematics, Nagoya University, 
Chikusa-ku, Furo-cho, Nagoya, 464-8602, Japan }
\email{furusho@math.nagoya-u.ac.jp}

\newtheorem{thm}{Theorem}[section]
\newtheorem{lem}[thm]{Lemma}

\newtheorem{prop}[thm]{Proposition}  

\theoremstyle{remark}
\newtheorem{ack}{Acknowledgments}        

\theoremstyle{definition}
\newtheorem{defn}[thm]{Definition}
\newtheorem{rem}[thm]{Remark}
\newtheorem{nota}[thm]{Notation}     
\newtheorem{note}[thm]{Note}
\newtheorem{eg}[thm]{Examples}       
\newtheorem{pf}{Proof}

\newtheorem{prob}[thm]{Problem}    
\newtheorem{assump}[thm]{Assumption}

\numberwithin{equation}{section}

\begin{document}
\bibliographystyle{amsalpha+}
\maketitle

\begin{abstract}
We establish a tannakian formalism of $p$-adic multiple polylogarithms
and $p$-adic multiple zeta values introduced in our previous paper
via a comparison isomorphism between 
a de Rham fundamental torsor and a rigid fundamental torsor 
of the projective line minus three points
and also discuss its Hodge and \'{e}tale analogues.
As an application we give a way to erase log poles of 
$p$-adic multiple polylogarithms and introduce
overconvergent $p$-adic multiple polylogarithms
which might be $p$-adic multiple analogue of 
Zagier's single-valued complex polylogarithms. 
\end{abstract}

\tableofcontents

\setcounter{section}{-1}
\section{Introduction}
This paper is the continuation of our previous paper \cite{F1}.
Let $p$ be a prime number. Let $m,k_1,\cdots,k_m\in\bold N$.
In \cite{F1} we constructed $p$-adic multiple polylogarithm
$Li_{k_1,\cdots,k_m}(z)$ ($z\in\bold C_p$), that is, 
a $p$-adic analogue of the (one-variable) complex multiple polylogarithm
defined locally near $0$ by
\begin{equation}\label{purple}
Li_{k_1,\cdots,k_m}(z)=\sum_{0<n_1<\cdots<n_m}
\frac{z^{n_m}}{n_1^{k_1}\cdots n_m^{k_m}}
\end{equation}
and $p$-adic multiple zeta value $\zeta_p(k_1,\cdots,k_m)$, 
that is, a $p$-adic analogue of multiple zeta value
\begin{equation}\label{gold}
\zeta(k_1,\cdots,k_m)=\sum_{0<n_1<\cdots<n_m}
\frac{1}{n_1^{k_1}\cdots n_m^{k_m}}\qquad (k_m>1).
\end{equation}
This was achieved by giving a $p$-adic Drinfel'd associator
$\Phi^p_\mathrm{ KZ}(A,B)$ of the $p$-adic KZ equation
and its fundamental solution $G_0(A,B)(z)$ ($z\in\bold C_p$).
The purpose of this paper is to relate these constructions to the motivic
fundamental torsors of the projective line minus three points.
To do this we will give a tannakian interpretation of the constructions
of \cite{F1}.
Using the Frobenius action on the rigid realization of the above torsor,
we will define the overconvergent $p$-adic multiple polylogarithm.
This function can be expressed by the $p$-adic multiple polylogarithm
constructed in \cite{F1}.
Using Hodge and $l$-adic \'{e}tale realizations we will recover the
complex multiple polylogarithm and the $l$-adic polylogarithms of
Wojtkowiak.\par

Our main tool is a comparison isomorphism (Lemma \ref{Austria}) 
\begin{equation}\label{clock}
\pi_1^{\mathrm{DR}}(\bold P^1_{\bold Q_p}\backslash\{0,1,\infty\}
:\overrightarrow{01},z)
\cong
\pi_1^{p,\mathrm{rig}}(\bold P^1_{\bold F_p}\backslash\{0,1,\infty\}
:\overrightarrow{01},z_0)
\end{equation}
between the de Rham (\S \ref{de Rham}) and the rigid (\S \ref{crys})
fundamental torsor where 
$z_0\in\bold P^1({\bold F_p})\backslash\{0,1,\infty\}$
is the reduction of
$z\in\bold P^1({\bold Q_p})\backslash\{0,1,\infty\}$ and
$\overrightarrow{01}$ is a tangential basepoint.
There is a canonical de Rham path $d_z$ (Lemma \ref{Akita}) in LHS
constructed from the canonical extension of unipotent connections
and a canonical rigid path $c_{z_0}$ (Lemma \ref{Niigata}) in RHS, 
Besser-Vologodsky's Frobenius invariant path \cite{Bes,Vo}.
Identifying these two torsors by \eqref{clock} 
we get a de Rham loop $d_z^{-1}c_{z_0}\in
\pi_1^{\mathrm{DR}}(\bold P^1_{\bold Q_p}\backslash\{0,1,\infty\}
:\overrightarrow{01})(\bold Q_p)$.
Our tannakian interpretation of $p$-adic multiple polylogarithms is

{\bf Theorem \ref{Sweden}.}
{\it
By the embedding 
$i:\pi_1^{\mathrm{DR}}(\bold P^1_{\bold Q_p}\backslash\{0,1,\infty\}
:\overrightarrow{01})(\bold Q_p)\hookrightarrow
\bold Q_p\langle\langle A,B\rangle\rangle$ \eqref{navy},
the loop $d_z^{-1}c_{z_0}$ corresponds to 
the non-commutative formal power series $G_0(A,B)(z)$ 
of the $p$-adic multiple polylogarithms.
}\\
Similarly by letting $z$ to be a tangential basepoint 
$\overrightarrow{10}$ we get
a tannakian interpretation of $p$-adic multiple zeta values:

{\bf Theorem \ref{Finland}.}
{\it
By the embedding $i$ 
the loop $d_{\overrightarrow{10}}^{-1}c_{\overrightarrow{10}}$
corresponds to 
the non-commutative formal power series $\Phi^p_\mathrm{ KZ}(A,B)$
of the $p$-adic multiple zeta values.
}\\
RHS of \eqref{clock} admits the Frobenius action $\phi_p$ 
(Definition \ref{Kanagawa}).
Deligne \cite{De2} introduced a variant of $p$-adic multiple zeta values 
to be a coefficient of the series
$\Phi^p_\mathrm{ De}(A,B)=
i(d_{\overrightarrow{10}}^{-1}\phi_p(d_{\overrightarrow{10}}))$
which we call the $p$-adic Deligne associator.

{\bf Theorem \ref{Netherland}.}
{\it
There is the following explicit relationship between 
the $p$-adic Drinfel'd associator and the $p$-adic Deligne associator:
$$
\Phi^p_{\mathrm{KZ}}(A,B)= \Phi^p_{\mathrm{De}}(A,B)\cdot
\Phi^p_{\mathrm{KZ}}\left(\frac{A}{p},
\Phi^p_{\mathrm{De}}(A,B)^{-1}\frac{B}{p}\Phi^p_{\mathrm{De}}(A,B)\right).  
$$
}\\
This gives formulae expressing our $p$-adic multiple zeta values 
in terms of Deligne's $p$-adic multiple zeta values 
and vice versa.
Namely
his $p$-adic multiple zeta values are equivalent to ours.
We introduce a pro-unipotent overconvergent iso-crystal
$({\mathcal V}^\dag_{\mathrm{KZ}},\nabla^\dag_{\mathrm{KZ}})$
on $\bold P^1_{\bold F_p}\backslash\{0,1,\infty\}$
associated with the $p$-adic KZ equation
and show

{\bf Proposition \ref{Croatia}.}
{\it
The pro-object 
$({\mathcal V}^\dag_{\mathrm{KZ}},\nabla^\dag_{\mathrm{KZ}})$
naturally admits a Frobenius structure.
}\\
The canonical section $1$ of ${\mathcal V}^\dag$ under 
the Frobenius action will map to a non-commutative power series 
$G^\dag_0(z)$ whose coefficient will be denoted by 
$Li^\dag_{k_1,\cdots,k_m}(z)$.
This series satisfies the differential equation in \cite{U,Y}
$$
dg=\left(\frac{A}{z}+\frac{B}{z-1}\right)gdz
-g\left(\frac{dz^p}{z^p}\frac{A}{p}+
\frac{dz^p}{z^p-1}
\Phi^p_{\mathrm{De}}(A,B)^{-1}\frac{B}{p}\Phi^p_{\mathrm{De}}(A,B)
\right)
$$
and $Li^\dag_{k_1,\cdots,k_m}(z)$ are an overconvergent analogue of our  
$p$-adic multiple polylogarithms.
The special value of these functions at $z=1$ will be
the Deligne's $p$-adic multiple zeta values. 
The relation between $G_0^\dag(A,B)(z)$ and our $G_0(A,B)(z)$
is given by the following:

{\bf Theorem \ref{Czech}.}
{\it
The overconvergent $p$-adic multiple polylogarithm 
$Li^\dag_{k_1,\cdots,k_m}(z)$ is expressed as a combination of our
$p$-adic multiple zeta values, $p$-adic multiple polylogarithms and
the $p$-adic logarithm by
$$
G_0^\dag(A,B)(z)=
G_0(A,B)(z)\cdot
G_0\left(\frac{A}{p} \ ,
\Phi^p_{\mathrm{De}}(A,B)^{-1}\frac{B}{p}\Phi^p_{\mathrm{De}}(A,B)\right)(z^p)^{-1}.  
$$
}\\
The function $Li_{k_1,\cdots,k_m}(z)$ is not overconvergent
but a Coleman function with log poles around $z=1$ and $\infty$.
The above formulae 
is a way to erase these log poles.
We will give calculations in Example \ref{Spain}. \par

The organization of this paper is as follows.
\S\ref{realizations} is preliminary for the rest of paper.
We will review the definition of various 
(de Rham, rigid, Betti and \'{e}tale) fundamental groups, torsors
with their various additional structures, i.e. 
the Frobenius action, the infinity Frobenius action and 
the Galois group action.
We also discuss tangential basepoints.
The Deligne's canonical path in the de Rham fundamental
torsor in Lemma \ref{Akita}
and Besser-Vologodsky's Frobenius invariant path 
in the rigid fundamental torsor in Lemma \ref{Niigata}
will be recalled.\par

\S \ref{results} is devoted to tannakian formalisms
which contains our main results.
In Theorem \ref{Sweden} and Theorem \ref{Finland}
we will clarify tannakian origins of the
$p$-adic multiple polylogarithms and $p$-adic multiple zeta values. 
An explicit relationship between our 
$p$-adic multiple zeta values and 
Deligne's $p$-adic multiple zeta values \cite{De2} 
will be stated in Theorem \ref{Netherland}.
A formula to express the overconvergent 
$p$-adic multiple polylogarithms in terms of 
our $p$-adic multiple polylogarithms 
will be stated in Theorem \ref{Czech}.
In \S \ref{Hodge} the Hodge analogue of these results will be discussed.
The tannakian origin of complex
multiple polylogarithms will be discussed
in Proposition \ref{Azerbaijan} 
and for multiple zeta values in Proposition \ref{Georgia}.
The infinity Frobenius action on 
the Betti fundamental torsor of 
$\bold P^1({\bold C})\backslash\{0,1,\infty\}$ will be used to
introduce another variant of multiple polylogarithms 
\eqref{Cambridge}.
In Theorem \ref{Russia} it will be shown that these functions are
single-valued and real-analytic. 
In Proposition \ref{Moldova} we
show a formula analogous to Theorem \ref{Czech} to express
them in terms of complex multiple polylogarithms.
\S \ref{Artin} is a brief explanation 
of analogous story in the \'{e}tale side.
We discuss an \'{e}tale analogue of polylogarithm \eqref{darkgray}
expressed in terms of Wojtkowiak's $l$-adic polylogarithm 
and an \'{e}tale analogue of Riemann zeta values (Example \ref{Ecuador})
expressed in term of Soul\'{e} characters and 
cyclotomic characters following \cite{NW,I90}.\par

\S \ref{motive} is motivic.
The algebra generated by $p$-adic multiple zeta values
will be related to 
Drinfel'd's \cite{Dr} pro-algebraic group $\underline{GRT}_1$,
Racinet's \cite{R} pro-algebraic group $\underline{DMR}_0$
and the motivic Galois group $\pi_1(\mathcal{MT}(\bold Z))$.
In \S \ref{Drinfeld} we recall Drinfel'd's 
pro-algebraic bi-torsor of Grothendieck-Teichm\"{u}ller.
In \S \ref{Racinet}, we recall Racinet's \cite{R} pro-algebraic
bi-torsor made by double shuffle relations and briefly explain 
a partial analogous story to \S \ref{Drinfeld}.
\S \ref{Goncharov} is a review of the formalism of mixed Tate motives in \cite{DG}.
The torsor of motivic Galois group is related with the torsors in 
\S \ref{Drinfeld} and \S \ref{Racinet}.
We give some motivic interpretations on Zagier's conjecture on multiple zeta values
and Ihara's conjecture on Galois image.
We also explain a motivic way of proving double shuffle relations for 
$p$-adic multiple zeta values using \cite{DG} and \cite{Y}.

\begin{ack}
This work was supported by the NSF grants DMS-0111298.
The author would like to express particular thanks to 
Prof. Deligne for many suggestions.
He is also grateful to Amir Jafari and Go Yamashita 
for many useful discussion.
He is also grateful to the referee for comments.
\end{ack}

\section{Preliminaries}\label{realizations}
This section is preliminaries for the next section.
We give materials which are indispensable to our main results by
recalling the definitions of various 
(de Rham, rigid, Betti and \'{e}tale) fundamental groups and torsors
with tangential base points and 
discussing additional structures there.
\par
In this section, we concentrate on a curve 
$X_{K}=\overline{X}_{K}-D_{K}$,
where $\overline{X}_{K}$ is
a proper smooth, geometrically connected
curve $\overline{X}_{K}$ over a field $K$ and
$D_{K}$ is its divisor over $K$
(sometimes we omit $K$).

\subsection{de Rham setting}\label{de Rham}
In this subsection,
we will recall the definitions of the de Rham fundamental
group, torsor, path space,
the tangential base point and the canonical 
base point and see how we obtain the canonical de Rham path
in the de Rham fundamental torsor by the canonical base point,
all of which are developed in \cite{De}\S 12.

\begin{nota}\label{Hokkaido}
In this subsection we assume that
$K$ is a field of characteristic $0$.
We denote the category of the nilpotent part 
\footnote{
It means the full subcategory consisting of objects
which are iterated extensions of the unit objects.}
of the category of the
pair $(\mathcal{V},\nabla)$ of a coherent
$\mathcal{O}_X$-module sheaf $\mathcal{V}$ on $X_K$
and an integrable (i.e. $\nabla\circ\nabla=0$) connection 
$\nabla:\mathcal{V}\to\mathcal{V}\otimes\Omega^1_{X_K}$
by $\mathcal{NC}^{\mathrm{DR}}(X_K)$,
It forms a neutral tannakian category 
(for the basics, consult \cite{DM})
by \cite{De}.
For our quick review of these materials, see \cite{S1}.
\end{nota}

In the de Rham realization, we consider three fiber functors.
Suppose that $x$ is a $K$-valued point of $X$,
i.e. $x:Spec K\to X_{K}$.
The first one is $\omega_x:\mathcal{NC}^{\mathrm{DR}}(X_K)\to Vec_K$
($Vec_K$: the category of finite dimensional $K$-vector spaces)
which is the pull-back of $\mathcal{NC}^{\mathrm{DR}}(X_K)$ by $x$.
Actually this forms a fiber functor (consult \cite{DM} for its definition) 
by \cite{S1}\S 3.1.
The second fiber functor
is the one introduced in \cite{De}\S 15
(and developed to higher dimensional case in \cite{BF}):
Let  $s\in D$ and
put $T_s=T_s\overline{X}_{K}$ and 
$T_s^\times=T_s\backslash\{0\}$($\cong \bold G_m$),
where $T_s\overline{X}_{K}$ is the tangent vector space of 
$\overline{X}_{K}$ at $s$.
By the construction of the extension of unipotent integrable 
connections to tangent bundles in \cite{De}\S15.28-15.36,
we have a morphism of tannakian categories (see also \cite{BF})
\begin{equation}\label{Heisei}
Res_s:{\mathcal{NC}}^{\mathrm{DR}}(X_K)\to
{\mathcal{NC}}^{\mathrm{DR}}(T^\times_s)
\end{equation}
which we call the tangential morphism.
By pulling back, each $K$-valued point $t_s$ of $T^\times_s$ 
determines a fiber functor
$\omega_{t_s}:{\mathcal{NC}}^{\mathrm{DR}}(X_K)\to Vec_K$
which we call the {\sf tangential base point}.
The last fiber functor is the special one 
under the assumption
\begin{equation}\label{Showa}
H^1(\overline{X},{\mathcal O}_{\overline{X}})=0.
\end{equation}
(i.e. $\overline{X}$ is a curve with genus $0$):
Let
$\omega_\Gamma:{\mathcal{NC}}^{\mathrm{DR}}(X_K)\to Vec_K$
be a functor sending each object 
$(\mathcal{V},\nabla)\in {\mathcal{NC}}^{\mathrm{DR}}(X_K)$
to $\Gamma(\overline{X}_K,{\mathcal V}_{\mathrm{can}})$,
the global section of its canonical extension
$(\mathcal{V}_{\mathrm{can}},\nabla_{\mathrm{can}})$ 
into $\overline{X}_K$ \cite{De} \S 12.2.
This functor forms a fiber functor (\cite{De} \S 12.4)
because $\mathcal{V}_{\mathrm{can}}$ forms a trivial bundle over 
$\overline{X}_K$ due to the condition \eqref{Showa}
by \cite{De} Proposition 12.3.
We call $\omega_\Gamma$
the {\sf canonical base point}.
We note that the assumption \eqref{Showa} is necessary.
The functor $\omega_\Gamma$ would be no longer a fiber functor 
if $\overline{X}$ was a proper smooth curve with genus $g>0$.

\begin{lem}[\cite{De}\S15.52]\label{Akita}
Let $\omega_*$ and $\omega_{*'}$ be any fiber functors above.
Then under the assumption \eqref{Showa}
there is a canonical isomorphism 
$d_{**'}:\omega_*\to\omega_{*'}$.
\end{lem}

\begin{pf}
It is because $\mathcal{V}_{\mathrm{can}}$ forms a trivial bundle
for $(\mathcal{V},\nabla)\in\mathcal{NC}^{\mathrm{DR}}(X_K)$
that we get a canonical isomorphism
$$\omega_x(\mathcal{V},\nabla)=\mathcal{V}_{(x)}=\mathcal{V}_{{\mathrm{can}},(x)}\cong
{\mathcal O}_{{\overline X}_K,(x)}\underset{K}{\otimes}
\Gamma(\overline{X}_K,\mathcal{V}_{{\mathrm{can}}})\cong
{\mathcal O}_{{\overline X}_K,(x)}\underset{K}{\otimes}
\omega_\Gamma(\mathcal{V},\nabla)\cong
\omega_\Gamma(\mathcal{V},\nabla).$$
and similarly 
$$\omega_{t_s}(\mathcal{V},\nabla)
\cong
{\mathcal O}_{T^\times_s,(t_s)}\underset{K}{\otimes}
\Gamma(\overline{X}_K,\mathcal{V}_{{\mathrm{can}}})
\cong
\omega_\Gamma(\mathcal{V},\nabla).$$
Here $_{(x)}$ stands for the fiber of sheaves at $x$.
By factoring through $\omega_\Gamma$,
we obtain a canonical isomorphism between any two of them.
\qed
\end{pf}

The fiber functors $\omega_x$ and $\omega_{t_s}$ 
(and $\omega_\Gamma$ under the assumption \eqref{Showa})
play a role of `base points'.

\begin{defn}[\cite{De}]\label{Yamagata}
Let $\omega_*$ and $\omega_{*'}$ be any fiber functors above.
The {\sf de Rham fundamental group}
$\pi_1^{\mathrm{DR}}(X_K:*)$ is $\underline{Aut}^\otimes(\omega_*)$
and the {\sf de Rham fundamental torsor}
$\pi_1^{\mathrm{DR}}(X_K:*,*')$ is $\underline{Isom}^\otimes(\omega_*,\omega_{*'})$
(for $\underline{Aut}^\otimes$ and 
$\underline{Isom}^\otimes$, see \cite{DM}).
\end{defn}

The former is a pro-algebraic group over $K$ and
the latter is a pro-algebraic bitorsor over $K$
where $\pi_1^{\mathrm{DR}}(X_K:*)$ acts from the right and
$\pi_1^{\mathrm{DR}}(X_K:*')$ acts from the left
\footnote{
For any two loops $\gamma$ and $\gamma'$,
the symbol $\gamma'\cdot\gamma$ stands for the composition of two 
which takes $\gamma$ at first and takes $\gamma'$ next.
Throughout this paper, we keep this direction.
}
and by definition
$\pi_1^{\mathrm{DR}}(X_K:*,*)=\pi_1^{\mathrm{DR}}(X_K:*)$.
Particularly under the assumption \eqref{Showa},
we have a {\sf canonical de Rham path} 
$d_{*,*'}\in\pi_1^{\mathrm{DR}}(X_K:*,*')(K)$
by Lemma \ref{Akita}.
This path is compatible with the composition:
$d_{*',*''}\cdot d_{*,*'}=d_{*,*''}$.
We also note that
for any field extension $K'/K$
the category $\mathcal{NC}^{\mathrm{DR}}(X_{K'})$ is equivalent to
$\mathcal{NC}^{\mathrm{DR}}(X_K)\otimes_K K'$
which is the category whose set of morphisms, $K'$-linear space,
is replaced by the extension of the set of morphisms of 
$\mathcal{NC}^{\mathrm{DR}}(X_K)$,
$K$-linear space, by $K'$ (\cite{De}\S 10.41),
whence
$\pi_1^{\mathrm{DR}}(X_{K'}:*)\cong 
\pi_1^{\mathrm{DR}}(X_K:*)\times_K K'$
and
$\pi_1^{\mathrm{DR}}(X_{K'}:*,*')\cong
\pi_1^{\mathrm{DR}}(X_K:*,*')\times_K K'$.

To establish the notion of the de Rham path space, 
we recall the notion of affine group $\mathcal T$-schemes 
for a tannakian category $\mathcal T$
(\cite{De}\S 5.4 and \cite{DG}\S 2.6).
The category of {\sf affine $\mathcal T$-schemes} is the dual of the category
of ind-objects of $\mathcal T$ which are unitary commutative algebras, that is,
each object is an ind-object $A$ of $\mathcal T$
with a product $A\otimes A\to A$
and a unit $1\to A$ verifying usual axioms.
We denote the corresponding affine $\mathcal T$-schemes with $A$ by
$Spec A$.
An {\sf affine group $\mathcal T$-scheme} is a group object of the category of
affine $\mathcal T$-schemes, i.e. the spectrum of a commutative Hopf algebra.

\begin{defn}\label{Aizu}
The {\sf de Rham path space}
\footnote{In \cite{De} \S 6.13, it is called the fundamental groupoid.
}
$\mathcal{P}_{X_K}^{\mathrm{DR}}$ 
is an affine group
$\mathcal{NC}^{\mathrm{DR}}(X_K)\otimes
\mathcal{NC}^{\mathrm{DR}}(X_K)$
\footnote{
That is the tensor product of tannakian categories \cite{De}\S 5.18.
}
-scheme introduced in  \cite{De}\S 6.13 and \cite{DG}\S 4,
which satisfies
$$
(\omega\otimes\omega')(\mathcal{P}_{X_K}^{\mathrm{DR}})=
\underline{Isom}^\otimes(\omega,\omega')(K)
$$
for any fiber functors
$\omega,\omega':\mathcal{NC}^{\mathrm{DR}}(X_K)\to Vec_K$
and is represented by an ind-object
$\mathcal{A}^{\mathrm{DR}}$ of
$\mathcal{NC}^{\mathrm{DR}}(X_K)\otimes
\mathcal{NC}^{\mathrm{DR}}(X_K)$
i.e. $\mathcal{P}_{X_K}^{\mathrm{DR}}=Spec \mathcal{A}^{\mathrm{DR}}$.
The {\sf de Rham fundamental path space
$\mathcal{P}^{\mathrm{DR}}_{X_K,\omega}$ 
with the base point $\omega$}
is its pull-back by
$\omega\otimes id$,
which is represented by the ind-object
$\mathcal{A}_{\omega}^{\mathrm{DR}}:=
(\omega\otimes id)(\mathcal{A}^{\mathrm{DR}})$
of
$\mathcal{NC}^{\mathrm{DR}}(X_K)$.
\end{defn}

By definition, we have
\begin{equation}\label{Nichiro}
(\omega_*\otimes\omega_{*'})\mathcal{P}^{\mathrm{DR}}_{X_K}=
\omega_{*'}(\mathcal{P}^{\mathrm{DR}}_{X_K,\omega_*})=
\pi_1^{\mathrm{DR}}(X_K:*,*')(K).
\end{equation}

\subsection{rigid setting}\label{crys}
In this paper we take rigid fundamental groups (torsors) as
a crystalline realization of motivic fundamental groups (torsors).
We will recall the definition of the rigid fundamental 
group \cite{CLS,S2}, torsor, path space,
the tangential base point \cite{BF,De}
and the canonical Frobenius invariant path
\cite{Bes,Vo}.\par

\begin{nota}\label{Maebashi}
In this subsection we assume that $K$ is 
a non-archimedean local field of characteristic $0$.
We denote $V$ to be its valuation ring and
$k$ to be the residue field of characteristic $p>0$ with $q(=p^r)$-elements.
We also assume that $\overline{X}_{K}$ with $D_{K}$
is a generic fiber of a proper smooth geometrically connected
curve $\mathcal{X}$
with a divisor $\mathcal D$ over $V$.
We denote its special fiber over $k$ 
by $\overline{X}_0$ with $D_0$.
Put $X_0=\overline{X}_0\backslash D_0$ and its natural embedding
$j:X_0\hookrightarrow\overline{X}_0$.
We denote the category of the nilpotent part of the 
category of (overconvergent
\footnote{
We can omit it because it is shown in \cite{CLS}
that unipotent (nilpotent) isocrystals are overconvergent.
}
) isocrystals, which consists of
the pair $(\mathcal{V}^\dag,\nabla^\dag)$ 
of a coherent $j^\dag\mathcal{O}_{]\overline{X}_0[}$
\footnote{
We do not call $\dag$ \lq plus' but \lq dagger'. 
For the definition of $j^\dag\mathcal{O}_{]\overline{X}_0[}$,
see \cite{Bes}.
}
-module sheaf $\mathcal{V}^\dag$ and an integrable connection
$\nabla^\dag:\mathcal{V}^\dag\to
\mathcal{V}^\dag\otimes_{\mathcal{O}_{]\overline{X}_0[}}
\Omega^1_{]\overline{X}_0[}$,
by $\mathcal{NI}^\dag(X_0)$.
Here $]\overline{X}_0[$ stands for the tubular neighborhood \cite{Ber}
of $\overline{X}_0$. 
This category depends only on $X_0$ and
forms a neutral tannakian category 
and depends only on $X_0$ (cf. \cite{Ber,CLS,S2}).
For our quick review of these materials, see \cite{Bes}.
\end{nota}

In the rigid realization, we consider three fiber functors.
Suppose that $x_0$ is a $k$-valued point of $X_0$,
i.e. $x_0:Spec \ k\to X_0$.
The first one is $\omega_{x_0}:\mathcal{NI}^\dag(X_0)\to Vec_{K_0}$
which is the pull-back of $\mathcal{NI}^\dag(X_0)$ by $x_0$.
This forms a fiber functor and is described by 
$\omega_{x_0}({\mathcal V}^\dag,\nabla^\dag)=
{\mathcal V}^\dag(]x_0[)^{\nabla^\dag}=
\{v\in{\mathcal V}^\dag(]x_0[)\bigm|\nabla^\dag(v)=0\}$
for $({\mathcal V}^\dag,\nabla^\dag)\in\mathcal{NI}^\dag(X_0)$
on the tubular neighborhood of $]x_0[$ 
\cite{Bes}.
The second fiber functor is 
$\omega_{s_0}:\mathcal{NI}^\dag(X_0)\to Vec_K$ ($s_0\in D_0$)
introduced in \cite{BF} Definition 1.3:
This is described by
$\omega_{s_0}({\mathcal V}^\dag,\nabla^\dag)=
\{v\in{\mathcal V}^\dag(]s_0[)[\log^a z_{s_0}]\bigm|\nabla^\dag(v)=0\}$,
where $\log^a$ is a branch of $p$-adic logarithm
such that $\log^ap=a\in K$
and $z_{s_0}$ is a local parameter of $]s_0[$
(cf.\cite{F1} \S 2.1).
The last fiber functor is the rigid version of the tangential base point:
Put $T_{s_0}=T_{s_0}\overline{X}_0$ and
$T^\times_{s_0}=T_{s_0}\backslash\{0\}(\cong\bold G_m)$
where $T_{s_0}\overline{X}_0$
is the tangent space of $\overline{X}_0$ at $s_0$.
By the construction in \cite{BF}\S 2,
we have a morphism of tannakian categories
\begin{equation}\label{Meiji}
Res_{s_0}:\mathcal{NI}^\dag(X_0)\to\mathcal{NI}^\dag(T^\times_{s_0}),
\end{equation}
which we call the tangential morphism.
Thus each $k$-valued point $t_{s_0}$ of $T^\times_{s_0}$ determines a fiber functor 
$\omega_{t_{s_0}}:\mathcal{NI}^\dag(X_0)\to Vec_K$ 
by pulling back.
Shortly the construction of $Res_{s_0}$ developed in \cite{BF}
is the transmission of the morphism \eqref{Heisei}
into the rigid setting 
by the categorical equivalences 
${\mathcal{NC}}^{\mathrm{DR}}(X_{K})\cong
\mathcal{NI}^\dag(X_0)$
and
${\mathcal{NC}}^{\mathrm{DR}}(T^\times_s)\cong
\mathcal{NI}^\dag(T^\times_{s_0})$
in \cite{CLS,S2},
where $s_0$ is the reduction of $s$.

The fiber functors, $\omega_{x_0}$, $\omega_{s_0}$ and $\omega_{t_{s_0}}$,
play a role of `base points'.

\begin{defn}\label{Chiba}
Let $\omega_*$ and $\omega_{*'}$ be any fiber functors above.
The {\sf rigid fundamental group}
$\pi_1^{p,{\mathrm{rig}}}(X_0:*)$ is
$\underline{Aut}^\otimes(\omega_*)$
and the {\sf rigid fundamental torsor}
$\pi_1^{p,{\mathrm{rig}}}(X_0:*,*')$ is 
$\underline{Isom}^\otimes(\omega_*,\omega_{*'})$
(cf. \cite{CLS,S2}).
The {\sf rigid path space}
$\mathcal{P}^\dag_{X_0}$ 
is an affine group
$\mathcal{NI}^{\dag}(X_0)\otimes
\mathcal{NI}^{\dag}(X_0)$-scheme
which satisfies
$$
(\omega\otimes\omega')(\mathcal{P}_{X_0}^{\dag})=
\underline{Isom}^\otimes(\omega,\omega')(K)
$$
for any fiber functors
$\omega,\omega':\mathcal{NI}^{\dag}(X_0)\to Vec_K$
and is represented by an ind-object
$\mathcal{A}^\dag$ of
$\mathcal{NI}^{\dag}(X_0)\otimes
\mathcal{NI}^{\dag}(X_0)$,
i.e. $\mathcal{P}_{X_0}^{\dag}=Spec \mathcal{A}^\dag$.
The {\sf rigid path space
$\mathcal{P}^{\dag}_{X_0,\omega}$ with the base point $\omega$}
is its pull-back by
$\omega\otimes id$,
which is represented by the ind-object
$\mathcal{A}_\omega^\dag:=(\omega\otimes id)(\mathcal{A}^\dag)$
of $\mathcal{NI}^{\dag}(X_0)$.
\end{defn}

We note that $\pi_1^{p,{\mathrm{rig}}}(X_0:*,*')$ 
is a pro-algebraic bitorsor over $K$
where $\pi_1^{p,{\mathrm{rig}}}(X_0:*)$ acts from the right and
$\pi_1^{p,{\mathrm{rig}}}(X_0:*')$ acts from the left.
By definition we have 
$$
(\omega_*\otimes\omega_{*'})\mathcal{P}^{\dag}_{X_0}=
\omega_{*'}(\mathcal{P}^{\dag}_{X_0,\omega_*})=
\pi_1^{p,{\mathrm{rig}}}(X_0:*,*')(K).
$$\par
Let $F_q$ be $r$-th power of the absolute Frobenius automorphism 
on $X_0$ and
its induced automorphism on $T^\times_{s_0}$ $(s_0\in D_0)$.
The pull-back functors 
$F_q^*:\mathcal{NI}^{\dag}(X_0)\to\mathcal{NI}^{\dag}(X_0)$
and $F_q^*:\mathcal{NI}^{\dag}(T^\times_{s_0})
\to\mathcal{NI}^{\dag}(T^\times_{s_0})$
determine equivalences of tannakian categories.
We note that they are compatible with \eqref{Meiji} (\cite{BF} \S 2).
They induce isomorphisms (\cite{De} \S 10.44)
\begin{equation}\label{Edo}
F_{q}:{\mathcal P}_{X_0}\to F_q^*{\mathcal P}_{X_0},
\end{equation}
\begin{equation}\label{Azuchi}
F_{q}:{\mathcal P}_{X_0,\omega}\to 
F_q^*({\mathcal P}_{X_0,F_q^*\omega})
\end{equation}
\begin{equation}\label{Momoyama}
F_{q}:\pi_1^{p,{\mathrm{rig}}}(X_0:*,*')\to
\pi_1^{p,{\mathrm{rig}}}(X_0:*,*').
\end{equation}

\begin{defn}\label{Kanagawa}(\cite{De} \S 10.44)
The {\sf Frobenius actions} $\phi_q$ stands for the inverse
\footnote{
We take the inverse to follow \cite{De} \S 13.6.
}
$F_q^{-1}$ in \eqref{Edo} $\sim$ \eqref{Momoyama}.
\end{defn}

By \eqref{Edo} and \eqref{Azuchi}
the ind-objects $\mathcal{A}^\dag$ and $\mathcal{A}_\omega^{\dag}$
in Definition \ref{Chiba}
naturally admit structures of {\it ind}-overconvergent $F$-isocrystals,
where an {\sf overconvergent $F$-isocrystal} means an
overconvergent isocrystal $V$
with a {\sf Frobenius structure} $\phi$
(i.e. a horizontal isomorphism $\phi:F_q^*V\to V$).
The description of the Frobenius structure in case of
$X=\bold P^1\backslash\{0,1,\infty\}$
is one of our main topics in this paper.
The following lemma 
which takes a place of Lemma \ref{Akita} in our rigid setting
is indispensable to our study.

\begin{lem}\label{Niigata}
Let $\omega_*$ and $\omega_{*'}$ be any fiber functors defined above.
Then there exists a unique path 
$c_{*,*'}\in \pi_1^{p,{\mathrm{rig}}}(X_0:*,*')(K)$
which is invariant under the Frobenius action, i.e.
$F_q(c_{*,*'})=c_{*,*'}$.
\end{lem}

This result is shown by Besser \cite{Bes} and Vologodsky \cite{Vo}
for usual base points (more generally in higher dimensional setting)
but their proofs also work for tangential base points 
and points on $D_0$ directly.
We call $c_{*,*'}$ the {\sf Frobenius invariant path}.
We have compatibilities, such as functorialities and
$c_{*',*''}\cdot c_{*,*'}=c_{*,*''}$
(for more details, see \cite{Bes}).

\subsection{Betti setting}\label{Betti}
We will briefly recall the definition of the Betti fundamental group,
torsor and the tangential base point in \cite{De}\S15.

\begin{nota}\label{Yamanashi}
In this subsection, we assume $k=\bold C$.
We denote the category of the nilpotent part of the category 
of local systems
\footnote{
We mean covariant functors $V:\Pi_{X(\bold C)}\to Vec_{\bold Q}$
where $\Pi_{X(\bold C)}$ stands for the category
whose objects are points on $X(\bold C)$ and
whose set of morphisms between two points $*$ and $*'$ on 
$X(\bold C)$ is the usual topological fundamental torsor
$\pi_1^{\mathrm{top}}(X(\bold C):*,*')$.
}
of finite dimensional $\bold Q$-vector spaces 
over the topological space
$X(\bold C)$ by $\mathcal{NL}^{\mathrm{Be}}(X(\bold C))$.
It also forms a neutral tannakian category over $\bold Q$
by \cite{De}\S 10.
\end{nota}
 
In the Betti realization, we consider two fiber functors.
Suppose that $x\in X(\bold C)$.
The first one is 
$\omega_x:\mathcal{NL}^{\mathrm{Be}}(X(\bold C))\to Vec_{\bold Q}$
which associates each nilpotent local system $V$ 
with the pullback $V(x)$.
Actually this forms a fiber functor of 
$\mathcal{NL}^{\mathrm{Be}}(X(\bold C))$ by \cite{De}\S 10.
The second one is the tangential basepoint:
Let $s\in D(\bold C)$ and
put $T_s=T_s\overline{X}(\bold C)$ and 
$T_s^\times=T_s\backslash\{0\}$($\cong \bold G_m$),
where $T_s\overline{X}(\bold C)$ is the tangent vector space of 
$\overline{X}(\bold C)$ at $s$.
By \cite{De}\S\S 15.3-15.12,
we have a morphism of tannakian categories
\begin{equation}\label{Sengoku}
Res_s:\mathcal{NL}^{\mathrm{Be}}(X(\bold C))\to
\mathcal{NL}^{\mathrm{Be}}(T_s^\times(\bold C)).
\end{equation}
which we call the tangential morphism.
By pulling back, each point $t_s\in T_s^\times(\bold C)$
determines a fiber functor 
$\omega_{t_s}:\mathcal{NL}^{\mathrm{Be}}(X(\bold C))\to Vec_{\bold Q}$.

The above fiber functors, $\omega_x$ and $\omega_{t_s}$,
play a role of \lq base points' and
the following is a reformulation of the Betti fundamental group
and torsor in \cite{De}\S 13.

\begin{defn}\label{Nagano}
Let $\omega_*$ and $\omega_*'$ be any fiber functors above.
The {\sf Betti fundamental group}
$\pi_1^{\mathrm{Be}}(X(\bold C):*)$ is $\underline{Aut}^\otimes(\omega_*)$
and the {\sf Betti fundamental torsor}
$\pi_1^{\mathrm{Be}}(X(\bold C):*,*')$ is 
$\underline{Isom}^\otimes(\omega_*,\omega_{*'})$.
\end{defn}

The former is a pro-algebraic group over $\bold Q$ and
the latter is a pro-algebraic bi-torsor over $\bold Q$ 
where $\pi_1^{\mathrm{Be}}(X(\bold C):*)$ acts on the right 
and $\pi_1^{\mathrm{Be}}(X(\bold C):*')$ acts on the left.
We note that we have  
a natural morphism from the topological fundamental group
$\pi_1^{\mathrm{top}}(X(\bold C):*)$
(resp. torsor $\pi_1^{\mathrm{top}}(X(\bold C):*,*')$)
to the set of $\bold Q$-valued points on 
the Betti fundamental group $\pi_1^{\mathrm{Be}}(X(\bold C):*)$ 
(resp. torsor $\pi_1^{\mathrm{Be}}(X(\bold C):*,*')$),
which induces an isomorphism between 
$\pi_1^{\mathrm{Be}}(X(\bold C):*)$
and the Malcev (unipotent) completion \cite{De}\S 9
over $\bold Q$ of
$\pi_1^{\mathrm{top}}(X(\bold C):*)$ (\cite{De} \S 10.24).
As for the topological fundamental groups and torsors
with tangential base points,
see \cite{De}\S\S 15.3-15.12.
Although we have a canonical path, 
both in the de Rham setting $d_{*,*'}$ Lemma \ref{Akita} 
and in the rigid setting $c_{*,*'}$ Lemma \ref{Niigata},
the author do not know whether we have any canonical path
connecting any two base points $*$ and $*'$ 
in the Betti setting or not.

\begin{assump}\label{Aichi}
Further we assume that $(\overline{X}, D)$ has a real structure
$(\overline{X'}, D')$, i.e.
$(\overline{X}, D)=
(\overline{X'}\times_{\bold R}\bold C, D'\times_{\bold R}\bold C)$, 
\end{assump}

Let $F_\infty$ be an involution of $X(\bold C)$ induced from 
the complex conjugation $z\in\bold C\mapsto \bar{z}\in\bold C$
and its induced automorphism on $T^\times_s$ ($s\in D$).
The pull-back functor
$F^*_{\infty}:\mathcal{NL}^{\mathrm{Be}}(X(\bold C))\to
\mathcal{NL}^{\mathrm{Be}}(X(\bold C))$ and
$F^*_{\infty}:\mathcal{NL}^{\mathrm{Be}}(T^\times_s(\bold C))\to
\mathcal{NL}^{\mathrm{Be}}(T^\times_s(\bold C))$ 
determine equivalences of tannakian categories.
We note that they are compatible with \eqref{Sengoku}.
They induce the isomorphism 
$$
F_{\infty}:\pi_1^{\mathrm{Be}}(X(\bold C):*,*')\to
\pi_1^{\mathrm{Be}}(X(\bold C):\bar{*},\bar{*}')
$$
where $\bar{*}$ (resp. $\bar{*}'$) means 
$F_{\infty}(*)$ (resp. $F_{\infty}(*')$).

\begin{defn}\label{Nagoya}
The {\sf infinity Frobenius action} $\phi_{\infty}$
stands for the inverse of $F_{\infty}$ above.
\end{defn}

The description of the infinity Frobenius action 
in case of $X=\bold P^1\backslash\{0,1,\infty\}$
will be discussed in \S \ref{results}.

\subsection{\'{e}tale setting}\label{etale}
We will briefly recall the definition of 
the \'{e}tale fundamental group, torsor and
the tangential base point in \cite{De}\S 15.

\begin{nota}\label{Mie}
In this subsection, we assume that $K$ is a field of characteristic $0$.
We denote $\overline{K}$ to be its algebraic closure and
$X_{\overline{K}}$ to be the scalar extension of $X_K$ by $\overline{K}$.
Let $l$ be a prime 
(we do not need to assume that $l$ is inequal to 
the prime $p$ in \S \ref{crys}).
We denote the category of the nilpotent part of the category 
of \'{e}tale local systems
\footnote{
We mean covariant functors 
$V:\hat{\Pi}_{X_{\overline K}}\to Vec_{\bold Q_l}$
where $\hat{\Pi}_{X_{\overline K}}$ stands for the category
whose objects are points on $X(\overline K)$ and
whose set of morphisms between two points $*$ and $*'$ on 
$X(\overline K)$ is the pro-finite fundamental torsor
$\hat{\pi}_1(X_{\overline K}:*,*')$ in \cite{SGA1}.
}
of finite dimensional $\bold Q_l$-vector spaces 
over 
$X(\overline{K})$ by $\mathcal{NE}^l(X_{\overline K})$.
It also forms a neutral tannakian category over $\bold Q_l$.
\end{nota}

In the \'{e}tale realization, we consider two fiber functors.
Suppose that $x\in X(\overline{K})$.
The first one is the pullback
$\omega_x:\mathcal{NE}^l(X_{\overline K})\to Vec_{\bold Q_l}$
which associates each nilpotent \'{e}tale local system $V$
with the pullback $V(x)$.
The second one is the tangential base point:
Let $s\in D(\overline K)$ and 
put $T_s=T_s\overline{X}_{\overline K}$ and   
$T_s^\times=T_s\backslash\{0\}$($\cong \bold G_m$),
where $T_s\overline{X}_{\overline K}$ is the tangent vector space of 
$\overline{X}_{\overline K}$ at $s$.
By \cite{De}\S\S 15.13-15.26,
we have a morphism of tannakian categories
\begin{equation}\label{Muromachi}
Res_s:\mathcal{NE}^l(X_{\overline K})\to\mathcal{NE}^l(T^\times_s)
\end{equation}
which we call the tangential morphism.
By pulling back, each $\overline K$-valued point $t_s$ in $T^\times_s$ determines a fiber functor 
$\omega_{t_s}:\mathcal{NE}^l(X_{\overline K})\to Vec_{\bold Q_l}$.

The above fiber functors, $\omega_x$ and $\omega_{t_s}$,
play a role of \lq basepoints' and
the following is a reformulation of the $l$-adic \'{e}tale 
fundamental group and torsor in \cite{De} \S 13.

\begin{defn}\label{Shiga}
Let $\omega_*$ and $\omega_{*'}$ be any fiber functors above.
The {\sf $l$-adic \'{e}tale fundamental group}
$\pi_1^{l,{\text{\'{e}t}}}(X_{\overline K}:*)$ is $\underline{Aut}^\otimes(\omega_*)$
and the {\sf $l$-adic \'{e}tale fundamental torsor}
$\pi_1^{l,{\text{\'{e}t}}}(X_{\overline K}:*,*')$ is 
$\underline{Isom}^\otimes(\omega_*,\omega_{*'})$.
\end{defn}

The former is a pro-algebraic group over $\bold Q_l$ and
the latter is a pro-algebraic bitorsor over $\bold Q_l$
where $\pi_1^{l,{\text{\'{e}t}}}(X_{\overline K}:*)$ acts on the right and
$\pi_1^{l,{\text{\'{e}t}}}(X_{\overline K}:*')$ acts on the left.
Let $\overline{K'}$ be an algebraically closed field containing 
$\overline{K}$.
The scalar extension induces a categorical equivalence
$\mathcal{NE}^l(X_{\overline K'})\cong
\mathcal{NE}^l(X_{\overline K})$,
from which we get isomorphisms
$\pi_1^{l,{\text{\'{e}t}}}(X_{\overline K'}:*)\cong
\pi_1^{l,{\text{\'{e}t}}}(X_{\overline K}:*)$
and
$\pi_1^{l,{\text{\'{e}t}}}(X_{\overline K'}:*,*')\cong
\pi_1^{l,{\text{\'{e}t}}}(X_{\overline K}:*,*')$
(\cite{De} \S 10.19).\par

Suppose that $*$ and $*'$ are $K$-valued points on $X$
or $T^\times_s$ for $s\in D$
(this assumption is attained by enlarging the base field $K$).
Let $\sigma\in Gal(\overline{K}/K)$.
We get {\sf Galois action}
$$
\sigma:\pi_1^{l,{\text{\'{e}t}}}(X_{\overline K}:*,*')\to
\pi_1^{l,{\text{\'{e}t}}}(X_{\overline K}:*,*')
$$
from the Galois action 
$id\otimes\sigma:X_{\overline K}\to X_{\overline K}$
in a similar way to our construction of 
$F_{q*}$ in \S \ref{crys} and 
$F_{\infty *}$ in \S \ref{Betti}.
We note that
the pro-finite fundamental group $\hat{\pi_1}(X_{\overline K}:*)$ 
(resp. torsor $\hat{\pi_1}(X_{\overline K}:*,*')$)
also admits Galois action and we have a morphism  
from 
$\hat{\pi_1}(X_{\overline K}:*)\to
\pi_1^{l,\text{\'{e}t}}(X_{\overline K}:*)(\bold Q_l)$ 
(resp.
$\hat{\pi_1}(X_{\overline K}:*,*')\to
\pi_1^{l,\text{\'{e}t}}(X_{\overline K}:*,*')(\bold Q_l)$)
which is compatible with its group (resp. torsor) structure and 
their Galois group actions.
Actually this morphism induces an isomorphism between
$\pi_1^{l,\text{\'{e}t}}(X_{\overline K}:*)$
and the Malcev (unipotent) completion over $\bold Q_l$ of 
$\hat{\pi_1}(X_{\overline K}:*)$
\cite{De} \S 10.

\begin{lem}\label{Hyogo}
Assume that $K$ is a number field with ring of integers
${\mathcal O}_K$.
Let $*$ and $*'$ be any $K$-valued points.
Let $\frak p$ be a prime ideal, coprime to $l{\mathcal O}_K$,
where the triple $(X,*,*')$ has a proper smooth model 
over ${\mathcal O}_{K,\frak p}$ with a good reduction.
Then there exists a unique path
$e^\frak p_{*,*'}\in\pi_1^{l,\text{\'{e}t}}
(X_{\overline{K}}:*,*')(\bold Q_l)$
which is invariant under the action of Frobenius element
$Frob_{\frak p}$ of $\frak p$, i.e.
$Frob_{\frak p}(e^\frak p_{*,*'})=e^\frak p_{*,*'}$.
\end{lem}

\begin{pf}
By the above assumption
$\pi_1^{l,\text{\'{e}t}}(X_{\overline{K}}:*,*')$
is unramified at $\frak p$, whence
we have an action of the geometric Frobenius $Frob_{\frak p}$ 
of $\frak p$.
Since the eigenvalues of the $Frob_\frak p$-action
on $H^1_{\text{\'{e}t}}(X_{\overline{K}},\bold Q_l)$ 
are positive weight Weil numbers,
we obtain a unique Frobenius invariant path
$e^\frak p_{*,*'}$
by an inductive construction same to
\cite{Bes}\S 3 (or \cite{Vo}\S 4.3). 
\qed
\end{pf}

We note that this Frobenius invariant path $e^\frak p_{*,*'}$ equips 
nice compatibilities 
such as functorialities and 
$e^\frak p_{*',*''}\cdot e^\frak p_{*,*'}=e^\frak p_{*,*''}$
similar to $d_{*,*'}$ in \S \ref{de Rham} and
$c_{*,*'}$ in \S \ref{crys}.
As Besser \cite{Bes} gave
the Frobenius invariant path $c_{*,*'}$
an interpretation of Coleman's $p$-adic iterated integration
in the rigid setting,
the author thinks that there may be also a significance
on the above special path $e^\frak p_{*,*'}$ in the $l$-adic 
\'{e}tale setting.

\section{Tannakian interpretations}\label{results}
This section is divided into three sides.
The first side is the Berthelot-Ogus side (\S \ref{Berthelot-Ogus})
where we shall give our main results.
A tannakian interpretation of $p$-adic MPL's 
and $p$-adic MZV's ($p$: prime)
\footnote{
MPL stands for multiple polylogarithm and
MZV stands for multiple zeta values.
}
is given in Theorem \ref{Sweden} and Theorem \ref{Finland}.
By describing Frobenius action on the rigid fundamental group 
of the projective line minus three points
we deduce an explicit formula  
between our $p$-adic MZV's and Deligne's $p$-adic MZV's \cite{De2,DG}
in Theorem \ref{Netherland}.
Further by considering Frobenius action on a rigid path space
we introduce an overconvergent version of $p$-adic MPL 
in Definition \ref{Bosnia}.
Its expression in terms of our $p$-adic MPL's is given in
Theorem \ref{Czech}.
The second side is the Hodge side (\S \ref{Hodge}),
where results analogous to \S \ref{Berthelot-Ogus} are presented. 
A tannakian reinterpretation
of MPL's and MZV's are given in Proposition \ref{Azerbaijan}
and Proposition \ref{Georgia} following \cite{BD}.
Using the infinity Frobenius action on the 
Betti fundamental torsor of the projective line minus three points
we introduce a new version of MPL \eqref{Cambridge}. 
Its monodromy-free property is shown in Theorem \ref{Russia}.
Its description is given in Proposition \ref{Moldova}.
The third side is the Artin side (\S \ref{Artin}).
We discuss $l$-adic \'{e}tale analogues of 
MPL's and MZV's defined in a way that Wojtkowiak \cite{W}
introduced the $l$-adic polylogarithms.\par

Throughout this and next section,
we fix 
$X=\bold P^1\backslash\{0,1,\infty\}
=Spec\bold Q[t,\frac{1}{t},\frac{1}{t-1}]
\subset\overline{X}=\bold P^1$
and $D=\{0,1,\infty\}$.
We denote the local parameter of $T^\times_0$
induced from $t$ by $t_0$ and the point 
$t_0=1$ on $T^\times_0$ by $\overrightarrow{01}$.
Similarly we denote the local parameter of $T^\times_1$
induced from $1-t$ by $t_1$ and the point 
$t_1=1$ on $T^\times_1$ by $\overrightarrow{10}$.

\subsection{Berthelot-Ogus side --- main results}
\label{Berthelot-Ogus}
In this side, we shall use a Berthelot-Ogus-type comparison
isomorphism (Lemma \ref{Austria}) between
a de Rham fundamental torsor and a rigid fundamental torsor for $X$.
Our first results are a tannakian interpretation of $p$-adic MPL's
and $p$-adic MZV's in Theorem \ref{Sweden} and
Theorem \ref{Finland} respectively.
Secondly we consider Frobenius action 
on a rigid fundamental group.
The definition of Deligne's $p$-adic MZV's \cite{De2,DG}
which describes the Frobenius action is recalled.
By using a tangential morphism we give an explicit formula 
between his and ours in Theorem \ref{Netherland}.
Thirdly we consider Frobenius action on a rigid path space.
{\it Overconvergent $p$-adic MPL} is introduced in
Definition \ref{Bosnia}.
It describes the Frobenius action 
and takes Deligne's $p$-adic MZV as \lq a special value' at $z=1$.
A formula to express it in terms of our $p$-adic MPL's is given in 
Theorem \ref{Czech}.


Before we state our first results
we give notations and basic lemma below
which will be frequently employed in this subsection.

\begin{nota}\label{Denmark}
We denote the special path 
$d_{\overrightarrow{01},z}$ (\S \ref{de Rham}) in 
the de Rham fundamental torsor
$\pi_1^\mathrm{ DR}(X:\overrightarrow{01},z)(\bold Q_p)$ 
for $z\in X(\bold Q_p)$ by $d_z$.
In particular we simply denote 
$d_{\overrightarrow{01},\overrightarrow{10}}$ in
by $d$.
Put $X_0=\bold P^1_{\bold F_p}\backslash\{0,1,\infty\}\subset
\overline{X_0}=\bold P^1_{\bold F_p}$.
Similarly we denote the special path $c_{\overrightarrow{01},z_0}$ 
(\S \ref{crys})
in the rigid fundamental torsor
$\pi_1^{p,\mathrm{ rig}}(X_0:\overrightarrow{01},z_0)(\bold Q_p)$ 
for $z_0\in\overline{X_0}(\bold F_p)$ by $c_{z_0}$.
In particular we simply denote
$c_{\overrightarrow{01},\overrightarrow{10}}$ by $c$.
The composition of the embedding
$\pi_1^\mathrm{ DR}(T_0^\times:\overrightarrow{01})\to
\pi_1^\mathrm{ DR}(X:\overrightarrow{01})$ by \eqref{Heisei}
and the exponential map
$H_1^\mathrm{ DR}(T_0^\times,\bold Q)\to
\pi_1^\mathrm{ DR}(T_0^\times,\overrightarrow{01})(\bold Q)$
sends the dual of $\frac{dt_0}{t_0}$ to 
a loop denoted by $x$.
Similarly a composite morphism to 
$\pi_1^\mathrm{ DR}(T_0^\times,\overrightarrow{10})(\bold Q)$
sends the dual of $\frac{dt_1}{t_1}$ to 
a loop denoted by $y'$.
Put $y=d^{-1}y'd\in\pi_1^\mathrm{ DR}(X:\overrightarrow{01})(\bold Q)$ 
Then it is easy to see that $x$ and $y$ are free generators of the
pro-algebraic group $\pi_1^\mathrm{ DR}(X:\overrightarrow{01})$.
In this section, we consider the embedding into 
the non-commutative formal power series 
with two variables $A$ and $B$
\begin{equation}\label{navy}
i:\pi_1^\mathrm{DR}(X:\overrightarrow{01})(\bold Q)\hookrightarrow
\bold Q\langle\langle A,B \rangle\rangle
\end{equation}
which sends 
$x\mapsto e^A:=1+\frac{A}{1!}+\frac{A^2}{2!}+\cdots$ and
$y\mapsto e^B:=1+\frac{B}{1!}+\frac{B^2}{2!}+\cdots$.
We note that this morphism is isomorphic to the embedding of
$\pi_1^\mathrm{DR}(X:\overrightarrow{01})(\bold Q)$
into the universal embedding algebra of 
$Lie\pi_1^\mathrm{DR}(X:\overrightarrow{01})(\bold Q)$.
By abuse of notation
we also employ the same symbol $i$ when we mention its 
scalar extension. 
\end{nota}

The following Berthelot-Ogus type comparison isomorphism
for fundamental torsors is our basics in this subsection.

\begin{lem}\label{Austria}
Let $z$ be a $\bold Q_p$-valued point of $X$ 
or $T^\times_s$ (for $s\in D$)
and $z_0$ be its modulo $p$ reduction. 
\begin{enumerate}
\item
If $z_0$ lies on $X_0(\bold F_p)$ or $T^\times_{s_0}(\bold F_p)$,
\begin{equation*}
\pi_1^{\mathrm{DR}}(X_{\bold Q_p}:\overrightarrow{01},z)
\cong
\pi_1^{p,\mathrm{rig}}(X_0:\overrightarrow{01},z_0).
\end{equation*}
\item
If $z_0$ does not lie on $X_0(\bold F_p)$ nor 
$T^\times_{s_0}(\bold F_p)$,
\begin{equation*}
\pi_1^{\mathrm{DR}}(X_{\bold Q_p}:\overrightarrow{01},z)
\times\bold Q_{p,st}\cong
\pi_1^{p,\mathrm{rig}}(X_0:\overrightarrow{01},z_0)
\times\bold Q_{p,\mathrm{st}}
\end{equation*}
where $\bold Q_{p,\mathrm{st}}$
stands for the polynomial algebra $\bold Q_p[l(p)]$ with one variable $l(p)$.
\end{enumerate}
\end{lem}

\begin{pf}
An equivalence of tannakian categories
\begin{equation}\label{red}
\mathcal{NC}^\mathrm{DR}(X_{\bold Q_p})\cong
\mathcal{NI}^\dag(X_{\bold F_p})
\end{equation}
is shown in \cite{CLS,S2} in more general situation.
In case of  $z_0\in X_0(\bold F_p)$
it gives rise to a natural isomorphism
$$
\omega_z(\mathcal{V},\nabla)\cong
\omega_{z_0}(\mathcal{V}^\dag,\nabla^\dag)
$$
for $(\mathcal{V},\nabla)\in\mathcal{NC}^\mathrm{DR}(X_{\bold Q_p})$,
where $(\mathcal{V}^\dag,\nabla^\dag)\in\mathcal{NI}^\dag(X_{\bold F_p})$
is associated with $(\mathcal{V},\nabla)$ via \eqref{red}.
The equivalences \eqref{red} and
$\mathcal{NC}^\mathrm{DR}(T_{s, \bold Q_p}^\times)\cong
\mathcal{NI}^\dag(T_{s_0}^\times)$ are compatible with tangential morphisms
(cf. \cite{BF}).
They also give rise to an natural isomorphism
$$
\omega_z\Bigl(Res(\mathcal{V},\nabla)\Bigr)\cong
\omega_{z_0}\Bigl(Res(\mathcal{V}^\dag,\nabla^\dag)\Bigr)
$$
for $(\mathcal{V},\nabla)\in\mathcal{NC}^\mathrm{DR}(X_{\bold Q_p})$
in case of $z_0\in T_0^\times(\bold F_p)$.
The equivalence \eqref{red} and
the above two isomorphisms of fiber functors
follow the isomorphism in (1).
In the rest case a natural isomorphism
$$
\omega_z(\mathcal{V},\nabla)\otimes\bold Q_p[l(p)]\cong
\omega_{z_0}(\mathcal{V}^\dag,\nabla^\dag)\otimes\bold Q_p[l(p)]
$$
is explained in \cite{Vo} \S 4.4, 
which follows the isomorphism in (2).
\qed
\end{pf}

We note that
this variable $l(p)$ reflects a branch parameter $a\in\bold Q_p$
of the $p$-adic polylogarithm $\log^a$.
By identifying the above two fundamental torsors by Lemma \ref{Austria}
we get a de Rham loop $d_z^{-1}c_{z_0}$,
which lies in 
$\pi_1^\mathrm{DR}(X:\overline{01})(\bold Q_p)$
in case of (1) and
in $\pi_1^\mathrm{DR}(X:\overline{01})(\bold Q_{p,\mathrm{ st}})$
in case of (2).

\begin{thm}\label{Sweden}
Let $z$ be a $\bold Q_p$-valued point of $X$ and
$z_0$ be its modulo $p$ reduction on $\overline{X_0}$. 
The de Rham loop $d_z^{-1}c_{z_0}$ 
corresponds to the special value of the fundamental solution $G_0$
\footnote{
In case of $z_0\in X_0(\bold F_p)$
it means the fundamental solution of
the $p$-adic KZ (Knizhnik-Zamolodchikov) equation
in $\bold Q_p\langle\langle A,B\rangle\rangle$ 
constructed in \cite{F1}Theorem 3.3
and
in case of $z_0\not\in X_0(\bold F_p)$
it means the series in 
$\bold Q_{p,\mathrm{ st}}\langle\langle A,B\rangle\rangle$ 
whose specialization at $l(p)=a$ is $G^a_0(z)$ in loc.cit.
}
at $z$ of the $p$-adic KZ equation
by the embedding $i$:
$$
i(d_{z}^{-1}c_{z_0})=G_0(z).
$$
\end{thm}

\begin{pf}
Our main tools to prove the theorem are 
a tannakian description of the KZ equation and 
Besser's tannakian interpretation of 
Coleman's $p$-adic iterated integration theory.\par

We start with the KZ equation.
Let $(\mathcal{V}_\mathrm{KZ},\nabla_\mathrm{KZ})$
be a pro-object of $\mathcal{NC}^\mathrm{DR}(X)$
associated with the KZ equation.
Namely  
${\mathcal V}_{\mathrm{KZ}}=
{\mathcal O}_X\langle\langle A,B\rangle\rangle$
and $\nabla_{\mathrm{KZ}}(g)=dg-(\frac{A}{t}+\frac{B}{t-1})g dt$
for local section $g$. 
By \cite{De} Proposition 12.10,
the Lie algebra $\frak p^{\mathrm{DR}}$ of $\pi_1^{\mathrm{DR}}(X:\Gamma)$
is free generated by two representatives 
$(\frac{dt}{t})^*$ and $(\frac{dt}{t-1})^*$
in $H_1^{\mathrm{DR}}(X,\bold Q)$
which are dual basis of $\frac{dt}{t}$ and $\frac{dt}{t-1}$ in 
$H^1_{\mathrm{DR}}(X,\bold Q)\cong
H^0(\overline{X},\Omega^1_{\overline X}(\log D))$ respectively
because $H^1(\overline{X},{\mathcal O}_{\overline X})=0$.
By the basic theorem \cite{DM} of tannakian categories,
the fiber functor $\omega_\Gamma$
induces an equivalence of categories
\begin{equation}\label{Taisho}
{\mathcal{NC}}^{\mathrm{DR}}(X_K)\cong
{\mathcal N}Rep \ \frak p^{\mathrm{DR}}
\end{equation}
where the right hand side stands for 
the category of nilpotent representations of 
the Lie algebra on finite dimensional $\bold Q$-vector spaces.

\begin{lem}\label{Italy}
By \eqref{Taisho}
the pro-object $({\mathcal V}_{\mathrm{KZ}},\nabla_{\mathrm{KZ}})$
corresponds to the vector space
$\bold Q\langle\langle A,B \rangle\rangle$
with the left multiplication 
$L:\frak p^{\mathrm{DR}}\to
{\frak gl}(\bold Q\langle\langle A,B \rangle\rangle)$
such that 
$L\left((\frac{dt}{t})^*\right)=L_A$ and 
$L\left((\frac{dt}{t-1})^*\right)=L_B$.
Here $L_F$ for $F\in\bold Q\langle\langle A,B \rangle\rangle$
means the left multiplication by $F$.
\end{lem}

\begin{pf}
In \cite{De}\S 12.5,
Deligne give a recipe of the correspondence \eqref{Taisho}
(in more general situation):
Let $(\mathcal{V},\nabla)\in {\mathcal{NC}}^{\mathrm{DR}}(X)$.
We put $V=\Gamma(\overline{X},{\mathcal V}_{\mathrm{can}})=
\omega_\Gamma(\mathcal{V},\nabla)$
and decompose as
$\nabla_{\mathrm{can}}=d+\omega$
where
$d:V\otimes\mathcal{O}_{\overline{X}}\to 
V\otimes\Omega^1_{\overline{X}}
(\log D)$
is the differential induced from 
$d:\mathcal{O}_{\overline{X}}\to \Omega^1_{\overline{X}}(\log D)$
and $\omega\in H^0(\overline{X},\Omega^1_{{\overline X}}(\log D))
\otimes End V$.
By the integrability of $\nabla$, 
we have the Lie algebra homomorphism
$\rho:{\frak p}^{\mathrm{DR}}
\to{\frak {gl}}V$ 
such that 
$\rho|_{H^{\mathrm{DR}}_1(X,\bold Q)}=-\omega$,
which gives an object 
$(V,\rho)\in\mathcal{N}Rep \ {\frak p}^{\mathrm{DR}}$.
Conversely let
$(V,\rho)\in\mathcal{N}Rep \ {\frak p}^{\mathrm{DR}}$.
We put $\omega=-\rho|_{H_1^{\mathrm{DR}}(X,\bold Q)}(\alpha)\in 
H^1_{\mathrm{DR}}(X,\bold Q)\otimes End V\cong
H^0(\overline{X},\Omega^1_{{\overline X}}(\log D))
\otimes End V$ 
where $\alpha\in 
H^1_{\mathrm{DR}}(X,\bold Q)\otimes H^{\mathrm{DR}}_1(X,\bold Q)$
is a canonical tensor representing the identity on
$H_1^{\mathrm{DR}}(X,\bold Q)$.
This gives 
$({\mathcal V},\nabla)=
(V\otimes{\mathcal O}_{X},d_{X}+\omega|_{X})
\in{\mathcal{NC}}^{\mathrm{DR}}(X)$.

Following his recipe we see that 
$L$ corresponds to 
$$
\omega=\frac{dt}{t}L_A+\frac{dt}{t-1}L_B\in
H_1^{\mathrm{DR}}(X,\bold Q)\widehat{\otimes}{\frak gl}
(\bold Q\langle\langle A,B \rangle\rangle)
$$
which gives the claim.
\qed
\end{pf}

By abuse of notation, we denote its induced group homomorphism by
\begin{equation*}
L:\pi_1^{\mathrm{DR}}(X:\Gamma)\to Aut\bold Q
\langle\langle A,B \rangle\rangle.
\end{equation*}
By the constructions of $i$ and $L$ above,
we have
\begin{equation}\label{brown}
i=ev_1\circ L\circ Int \ d_{\overrightarrow{01},\Gamma}.
\end{equation}
Here 
$Int \ d_{\overrightarrow{01},\Gamma}:
\pi_1^\mathrm{DR}(X:\overrightarrow{01})(\bold Q)\to
\pi_1^\mathrm{DR}(X:\Gamma)(\bold Q)$
is the automorphism sending 
$\gamma\mapsto d_{\overrightarrow{01},\Gamma}\cdot\gamma\cdot
d_{\overrightarrow{01},\Gamma}^{-1}$
and
$ev_1:Aut
\bold Q\langle\langle A,B \rangle\rangle\to
\bold Q\langle\langle A,B \rangle\rangle$ 
is the evaluation at $1$
sending $\sigma\mapsto\sigma(1)$.

Our second tool is
Besser's tannakian interpretation of the Coleman's
$p$-adic iterated integration theory \cite{Col} in \cite{Bes}:
He called a set 
$v=\{v_{x_0}\in\omega_{x_0}(\mathcal{V}^\dag,\nabla^\dag)|
{x_0\in\overline{X_0}(\overline{\bold F_p})}\}$
for $(\mathcal{V}^\dag,\nabla^\dag)\in\mathcal{NI}^\dag(X_0)$
a collection of {\sf analytic continuation along Frobenius}
of horizontal section of $(\mathcal{V}^\dag,\nabla^\dag)$
if $c_{x_0,x'_0}(v_{x_0})=v_{x'_0}$ for all
$x_0,x'_0\in\overline{X_0}({\overline{\bold F_p}})$
where $c_{x_0,x'_0}$ is a Frobenius invariant path 
in Lemma \ref{Niigata}
and showed that each Coleman function of $X(\bold C_p)$
(in the sense of \cite{Bes} \S 5) is expressed as a set
$\{\theta(v_{x_0})\}_{x_0}$
with a 
$j^\dag{\mathcal O}_{]\overline{X_0}[}\otimes\bold C_p$-module
homomorphism
$\theta:{\mathcal V}^\dag\otimes\bold C_p\to 
j^\dag{\mathcal O}_{]\overline{X_0}[}\otimes\bold C_p$
and a collection $v=\{v_{x_0}\}_{x_0}$ 
of analytic continuation along Frobenius of a pair
$(\mathcal{V}^\dag,\nabla^\dag)\in\mathcal{NI}^\dag(X_0)$.
In \cite{BF} we further extend Coleman functions to normal bundles
by making use of Frobenius invariant path to tangential basepoints.

In our previous paper \cite{F1} we consider the $p$-adic KZ equation
and construct its fundamental solution $G_0(t)$,
a two variable non-commutative formal power series
with Coleman function coefficients.
In his terminologies it is a collection
$\left\{G_0|_{]x_0[}\right\}_{x_0}$
of analytic continuation along Frobenius of horizontal section of
$(\mathcal{V}^\dag_\mathrm{ KZ},\nabla^\dag_\mathrm{ KZ})$
(the associated pro-object of $\mathcal{NI}^\dag(X_0)$ with 
$(\mathcal{V}_\mathrm{ KZ},\nabla_\mathrm{ KZ})$),
i.e.
$c_{x_0,y_0}(G_0|_{]x_0[})=G_0|_{]y_0[}$
for $x_0$ and $y_0$.
And since in our case 
the tangential morphism \eqref{Heisei} for $s=0$
sends $(\mathcal{V}_\mathrm{KZ},\nabla_\mathrm{KZ})$ to
$(\mathcal{O}_{T^\times_0}\langle\langle A,B\rangle\rangle,
d-A\frac{dt_0}{t_0})$,
$t_0^A=1+\frac{\log^{a}t_0}{1!}A
+\frac{(\log^{a}t_0)^2}{2!}A^2
+\frac{(\log^{a}t_0)^3}{3!}A^3+\cdots$
might be a collection of analytic continuation along Frobenius
of horizontal section of 
$Res_0(\mathcal{V}^\dag_\mathrm{KZ},\nabla^\dag_\mathrm{KZ})$.
By our tangential base point interpretation (\cite{BF} Proposition 2.11)
of the notion of the constant term (\cite{BF} Definition 2.1)
of horizontal sections, we see that 
$t_0^A$ is an analytic continuation along Frobenius of $G_0(t)$
to $T^\times_0$ because
both the constant term of $G_0(t)$ at $t=0$
and the constant term of $t_0^A$ at $t_0=0$ 
are equal to $1$.
Namely
$c_{\overrightarrow{01},y_0}(t_0^A|_{]\overrightarrow{01}[})
=G_0(t)|_{]y_0[}$.

Let $z$ be a $\bold Q_p$-valued point $X$ with reduction $z_0$.
We regard $c_{z_0}$ to be a de Rham path
via the identification in Lemma \ref{Italy}.
The isomorphism 
$d_{*,\Gamma}:\omega_*(\mathcal{V}_\mathrm{ KZ},\nabla_\mathrm{ KZ})
\to\omega_\Gamma(\mathcal{V}_\mathrm{ KZ},\nabla_\mathrm{ KZ})=
\bold Q_p\langle\langle A,B\rangle\rangle$
for $*=\overrightarrow{01}$ and $z$
let us rewrite
\begin{equation}\label{green}
d_{z,\Gamma}\cdot c_{z_0}\cdot d_{\Gamma,\overrightarrow{01}}(1)=G_0(z).
\end{equation}

By \eqref{brown}, \eqref{green} and
$d_{z,\Gamma}=d_{\overrightarrow{01},\Gamma}\cdot d_{z}^{-1}$,
we prove Theorem \ref{Sweden}.
\qed
\end{pf}

In our previous paper \cite{F1} Theorem 3.15
we showed that in each coefficient of $G_0(z)$
there appears $p$-adic MPL
$Li_{k_1,\cdots,k_m}(z)$,
which is a Coleman function admitting the expansion \eqref{purple}
on $|z|_p<1$.
In precise the coefficient of $A^{k_m-1}B\cdots A^{k_1-1}B$ in $G_0(z)$
is $(-1)^mLi_{k_1,\cdots,k_m}(z)$:
$$
G_0(z)=1+\log zA-Li_1(z)B+
\cdots+(-1)^mLi_{k_1,\cdots,k_m}(z)A^{k_m-1}B\cdots A^{k_1-1}B
+\cdots.
$$
By Theorem \ref{Sweden} we may say that the loop $d_z^{-1}c_{z_0}$
is a tannakian origin of $p$-adic MPL's.

In case when $z$ is also a tangential base point $\overrightarrow{01}$,
similarly we have

\begin{thm}\label{Finland}
The de Rham loop 
$d^{-1}c\in\pi_1^\mathrm{DR}(X:\overrightarrow{01})(\bold Q_p)$
corresponds to the $p$-adic Drinfel'd associator 
$\Phi^p_\mathrm{KZ}\in\bold Q_p\langle\langle A,B\rangle\rangle$
by the embedding $i$, i.e.
$$
i(d^{-1}c)=\Phi^p_\mathrm{KZ}.
$$
\end{thm}

\begin{pf}
The $p$-adic Drinfel'd associator $\Phi^p_\mathrm{KZ}$
is the series constructed in \cite{F1} Definition 3.12
which is equal to the limit value (\cite{F1} Lemma 3.27)
\begin{equation}\label{rose}
{\underset{\epsilon\to 0}{\lim}^\prime} 
\exp\left(-\log^{a}\epsilon\cdot B\right)\cdot 
G^{a}_0(1-\epsilon)
\in\bold Q_p\langle\langle A,B\rangle\rangle.
\end{equation}
In other word it is 
the constant term of $G^a_0(z)$ at $z=1$.
Then by a similar argument to the proof of Theorem \ref{Sweden}
we have 
$c_{\overrightarrow{01},\overrightarrow{10}}
(t_0^A|_{]\overrightarrow{01}[})
=\Phi^p_\mathrm{KZ}\cdot
t_1^B|_{]\overrightarrow{10}[}$,
whence
\begin{equation}\label{white}
d_{\overrightarrow{10},\Gamma}\cdot c\cdot d_{\Gamma,\overrightarrow{01}}(1)
=\Phi^p_\mathrm{KZ}.
\end{equation}
By \eqref{brown}, \eqref{white} and
$d_{\overrightarrow{10},\Gamma}=d_{\overrightarrow{01},\Gamma}\cdot d^{-1}$,
we prove Theorem \ref{Finland}.
\qed
\end{pf}

In our previous paper \cite{F1} Theorem 3.30 
we showed that in each coefficient of 
the $p$-adic Drinfel'd associator $\Phi^p_\mathrm{KZ}$ 
there appears $p$-adic MZV 
$\zeta_p(k_1,\cdots,k_m)$ ($k_m>1$)
introduced in \cite{F1},
which is a $p$-adic analogue of \eqref{gold}.
In precise the coefficient of $A^{k_m-1}B\cdots A^{k_1-1}B$ in 
$\Phi^p_\mathrm{KZ}$ 
is $(-1)^m\zeta_p(k_1,\cdots,k_m)$:
$$
\Phi^p_\mathrm{KZ}=1+\cdots+
(-1)^m\zeta_p(k_1,\cdots,k_m)A^{k_m-1}B\cdots A^{k_1-1}B+\cdots.
$$
By Theorem \ref{Finland} we may say that the loop $d^{-1}c$
is a tannakian origin of $p$-adic MZV's.

\begin{note}\label{Scottland}
The $p$-adic Drinfel'd associator $\Phi^p_\mathrm{KZ}$ is group-like \cite{F1},
that means 
$\Delta\Phi^p_\mathrm{KZ}=\Phi^p_\mathrm{KZ}\widehat\otimes\Phi^p_\mathrm{KZ}$ 
where $\Delta:\bold Q_p\langle\langle A,B\rangle\rangle\to
\bold Q_p\langle\langle A,B\rangle\rangle$
is the linear map induced from $\Delta(A)=A\otimes  1+1\otimes A$ and
$\Delta(B)=B\otimes 1+1\otimes B$.
Hence each coefficient of $\Phi^p_\mathrm{KZ}$ must satisfy (integral)
shuffle product formula (cf. \cite{F1}, \cite{FJ}).
We can recover and express general coefficients of $\Phi^p_\mathrm{KZ}$ 
in terms of $\zeta_p(k_1,\cdots,k_m)$'s ($k_m>1$) by the following method:
Let $W$ be a word, i.e. a monic and monomial element in 
$\bold Q_p\langle\langle A,B\rangle\rangle$.
Put $C_W$ to be its coefficient of $\Phi^p_\mathrm{KZ}$.
We have $C_A=C_B=0$.
For a convergent word $W$ written as $A^{k_m-1}B\cdots A^{k_1-1}B$ ($k_m>1$),
$C_W=(-1)^m\zeta_p(k_1,\dots,k_m)$.
For a divergent word $W'$ written as $B^rW$ ($r>0$) with $W$: convergent,
shuffle product formula gives 
$$
C_B^r\cdot C_W=r!C_{B^rW}+\text{other terms}.
$$
By an induction with respect to $r$, 
$C_{W'}$ is calculated and is expressed in terms of 
$C_W$'s with $W$: convergent.
For a word $W''=B^rWA^s$ ($s>0$) with $W$: convergent,
a similar induction argument let us able to calculate $C_{W''}$.
\end{note}

Our second result is on a description of a Frobenius action.
We consider a Frobenius action $F_p$ on 
$\pi_1^{\mathrm{DR}}(X:\overrightarrow{01},\overrightarrow{10})
(\bold Q_p)$
by transmitting that on 
$\pi_1^{p,\mathrm{rig}}(X_0:\overrightarrow{01},\overrightarrow{10})
(\bold Q_p)$ by Lemma \ref{Austria}.
It gives a new path $F_p(d)$ in 
$\pi_1^{\mathrm{DR}}(X:\overrightarrow{01},\overrightarrow{10})
(\bold Q_p)$
which is our central object to discuss here.
In Arizona Winter School 2002, 
Deligne \cite{De2} introduced another version of $p$-adic MZV 
which has a different tannakian origin from ours 
(see also \cite{DG} \S 5).

\begin{defn}\label{England}
The {\sf $p$-adic Deligne associator} 
$\Phi^p_{\mathrm{De}}$ is the series in
$\bold Q_p\langle\langle A,B\rangle\rangle$
which corresponds to
$d^{-1}\phi_p(d)$:
$$
i(d^{-1}\phi_p(d))=\Phi^p_{\mathrm{De}}.
$$
\end{defn}

To adapt the notations to ours above,
we denote $\zeta_p^{\mathrm{De}}(k_1,\cdots,k_m)$ ($k_m>1$)
to be the coefficient of $A^{k_m-1}B\cdots A^{k_1-1}B$ 
in $\Phi^p_{\mathrm{De}}$ multiplied by $(-1)^m$:
$$
\Phi^p_{\mathrm{De}}=1+\cdots+(-1)^m
\zeta_p^{\mathrm{De}}(k_1,\cdots,k_m)A^{k_m-1}B\cdots A^{k_1-1}B
+\cdots
$$
and call it {\sf Deligne's $p$-adic MZV}.
Since $\Phi^p_{\mathrm{De}}$ is group-like, we can 
recover and express general coefficients of $\Phi^p_{\mathrm{De}}$ 
in terms of $\zeta_p^{\mathrm{De}}(k_1,\cdots,k_m)$'s ($k_m>1$)
as Note \ref{Scottland}.
Deligne's $p$-adic MZV's come from the loop $d^{-1}\phi_p(d)$
while ours come from the loop $d^{-1}c$.
So his are different from ours.
But we can give a close relationship between them:

\begin{thm}\label{Netherland}
\begin{equation}\label{Amsterdam}
\Phi^p_{\mathrm{KZ}}(A,B)= \Phi^p_{\mathrm{De}}(A,B)\cdot
\Phi^p_{\mathrm{KZ}}\left(\frac{A}{p},
\Phi^p_{\mathrm{De}}(A,B)^{-1}\frac{B}{p}\Phi^p_{\mathrm{De}}(A,B)\right).  
\end{equation}
Here the last term means the series substituting
$(\frac{A}{p},
\Phi^p_{\mathrm{De}}(A,B)^{-1}\frac{B}{p}\Phi^p_{\mathrm{De}}(A,B))$
in $\Phi^p_{\mathrm{KZ}}(A,B)$.
\end{thm}

This formula looks nasty.
But according to the convention
\footnote{
Let $k$ be a field with characteristic $0$. Let $c_1,c_2\in k$ and
$g_1(A,B), g_2(A,B)\in k\langle\langle A,B\rangle\rangle$.
The product $(c_2,g_2)\circ(c_1,g_1)$ is defined by 
$\Bigl(c_1c_2, g_2(A,B)g_1
\bigl(\frac{A}{c_2},g_2(A,B)^{-1}\frac{B}{c_2}g_2(A,B)\bigr)
\Bigr)$.
}
in \cite{F0}
of the Grothendieck-Teichm\"{u}ller group \cite{Dr}
we rewrite simpler
$(p,\Phi^p_{\mathrm{KZ}})=(p,\Phi^p_{\mathrm{De}})\circ
(1,\Phi^p_{\mathrm{KZ}})$, i.e.
\begin{equation}\label{teal}
(p,\Phi^p_{\mathrm{De}})=
(p,\Phi^p_{\mathrm{KZ}})\circ
(1,\Phi^p_{\mathrm{KZ}})^{-1}. 
\end{equation}

\begin{pf}
We extends the Frobenius action $\phi_p$ on 
$\pi_1^\mathrm{ DR}(X,\overrightarrow{01})(\bold Q_p)$
to $\bold Q_p\langle\langle A,B \rangle\rangle$
by $i$.

\begin{lem}\label{Ireland}
$\phi_p(A)=\frac{A}{p}$, \ $\phi_p(B)=
{\Phi^p_{\mathrm{De}}}^{-1}\frac{B}{p}\Phi^p_{\mathrm{De}}$.
\end{lem}

\begin{pf}
By the compatibility of the Frobenius action on the tangential morphism
$\bold Q_p(1)=
\pi_1^{p,\mathrm{rig}}(T_0^\times:\overrightarrow{01})(\bold Q_p)
\to\pi_1^{p,\mathrm{rig}}(X_0:\overrightarrow{01})(\bold Q_p)$
we have $\phi_p(x)=x^{\frac{1}{p}}$.
Similarly we have $\phi_p(y')=y'^{\frac{1}{p}}$.
Since Frobenius action is compatible with the torsor structure
we have
$
\phi_p(y)=\phi_p(d^{-1}y'd)=
(d^{-1}\phi_p(d))^{-1}d^{-1}y'^{\frac{1}{p}}d(d^{-1}\phi_p(d)) 
=(d^{-1}\phi_p(d))^{-1}y^\frac{1}{p}(d^{-1}\phi_p(d)). 
$
\qed
\end{pf}

By $d^{-1}\phi_p(d)=d^{-1}c\phi_p(d^{-1}c)^{-1}$
we have
$$
\Phi^p_{\mathrm{De}}= \Phi^p_{\mathrm{KZ}}\cdot
\phi_p(\Phi^p_{\mathrm{KZ}})^{-1}.
$$
By Lemma \ref{Ireland} we get Theorem \ref{Netherland}.
\qed
\end{pf}

By expanding our explicit formula,
we can express our $p$-adic MZV's arising from Frobenius invariant path
in terms of Deligne's $p$-adic MZV's describing Frobenius action 
very explicitly and vice versa.
The following  are the easiest examples.

\begin{eg}\label{France}
\begin{enumerate}
\item $\zeta_p^\mathrm{ De}(k)=(1-\frac{1}{p^k})\zeta_p(k)$ \quad ($k>1$).
\item $\zeta_p^\mathrm{ De}(a,b)=(1-\frac{1}{p^{a+b}})\zeta_p(a,b)
-(\frac{1}{p^b}-\frac{1}{p^{a+b}})\zeta_p(a)\zeta_p(b)
-\sum\limits_{r=0}^{a-1}(-1)^r(\frac{1}{p^{a-r}}-\frac{1}{p^{a+b}})
\binom{b-1+r}{b-1}
\zeta_p(a-r)\zeta_p(b+r)
-(-1)^a\sum\limits_{s=0}^{b-1}(\frac{1}{p^{b-s}}-\frac{1}{p^{a+b}})
\binom{a-1+s}{a-1}
\zeta_p(a+s)\zeta_p(b-s)$ \quad ($b>1$).
\end{enumerate}
\end{eg}
We may say that Deligne's $p$-adic MZV's are 
not equal to but equivalent to ours.

Our third result is on Frobenius structures.
We give a Frobenius structure on 
$({\mathcal V}^\dag_{\mathrm{KZ}},\nabla^\dag_{\mathrm{KZ}})$.
We introduce and discuss an overconvergent  variant of the $p$-adic MPL
which describes this structure.

\begin{prop}\label{Croatia}
The pro-object 
$({\mathcal V}^\dag_{\mathrm{KZ}},\nabla^\dag_{\mathrm{KZ}})$
naturally admits a Frobenius structure
\begin{equation}\label{orange}
\phi:
F^*_p({\mathcal V}^\dag_{\mathrm{KZ}},\nabla^\dag_{\mathrm{KZ}})\to
({\mathcal V}^\dag_{\mathrm{KZ}},\nabla^\dag_{\mathrm{KZ}}).
\end{equation}
\end{prop}

\begin{pf}
We just derive it from the Frobenius structure on 
${\mathcal A}^\dag_{\overrightarrow{01}}$
(Definition \ref{Chiba}).
The following arguments in the de Rham setting 
would help our understandings.

\begin{lem}\label{Germany}
The pro-object $({\mathcal V}_{\mathrm{KZ}},\nabla_{\mathrm{KZ}})$
is isomorphic to the pro-object 
$({\mathcal A}^{\mathrm{DR}}_{\overrightarrow{01}})^\lor$,
the dual of the ind-object 
${\mathcal A}^{\mathrm{DR}}_{\overrightarrow{01}}$
(Definition \ref{Aizu}).
\end{lem}

\begin{pf}
By the categorical equivalence \eqref{Taisho},
$({\mathcal A}^{\mathrm{DR}}_\Gamma)^\lor$ corresponds to 
a natural representation $\rho_\Gamma$ of 
$Lie \ \pi_1^{\mathrm{DR}}(X:\Gamma)$
on the universal enveloping algebra $U\pi_1^{\mathrm{DR}}(X:\Gamma)$
of $Lie \ \pi_1^{\mathrm{DR}}(X:\Gamma)$
(we note that $U\pi_1^{\mathrm{DR}}(X:\Gamma)$
is isomorphic to the dual $H$ of the coordinate ring
$\omega_\Gamma({\mathcal A}^{\mathrm{DR}}_\Gamma)$ of 
$\pi_1^{\mathrm{DR}}(X:\Gamma)$ by \cite{A} Theorem 2.5.3 because
$Lie \ \pi_1^{\mathrm{DR}}(X:\Gamma)$ is the set of primitive elements
in $H$ by loc.cit.\S 4.3).
By \eqref{Nichiro} this representation $\rho_\Gamma$
is induced from the left action of $\pi_1^{\mathrm{DR}}(X:\Gamma)$
on $\pi_1^{\mathrm{DR}}(X:\Gamma,\Gamma)$,
thus $U\pi_1^{\mathrm{DR}}(X:\Gamma)$ is a free 
$U\pi_1^{\mathrm{DR}}(X:\Gamma)$-module of rank $1$ by $\rho_\Gamma$.
Therefore we say that the corresponding representations $\rho_\Gamma$ and
$L$ (Lemma \ref{Italy}) are equivalent, whence we have  
$({\mathcal V}_{\mathrm{KZ}},\nabla_{\mathrm{KZ}})\cong
({\mathcal A}^{\mathrm{DR}}_\Gamma)^\lor$.
On the other hand by the canonical isomorphism 
$d_{\overrightarrow{01},\Gamma}:
\omega_{\overrightarrow{01}}\to\omega_\Gamma$
it is easy to see
$({\mathcal A}^{\mathrm{DR}}_{\overrightarrow{01}})^\lor\cong
({\mathcal A}^{\mathrm{DR}}_\Gamma)^\lor$.
Whence we get the claim.
\qed
\end{pf}

By the categorical equivalence \eqref{red}
we deduce 
$({\mathcal V}^\dag_{\mathrm{KZ}},\nabla^\dag_{\mathrm{KZ}})\cong
({\mathcal A}^\dag_{\overrightarrow{01}})^\lor$
from Lemma \ref{Germany}.
Consequently we get the Frobenius structure on 
$({\mathcal V}^\dag_{\mathrm{KZ}},\nabla^\dag_{\mathrm{KZ}})$
by transmitting from
$\phi:{\mathcal A}^\dag_{\overrightarrow{01}}\to
F_p^*{\mathcal A}^\dag_{\overrightarrow{01}}$.
We get Proposition \ref{Croatia}.
\qed
\end{pf}

In the following section we introduce single-valued real-analytic 
MPL $Li_{k_1,\cdots,k_m}^-(z)$.
Here we introduce the following $p$-adic analogue
$Li_{k_1,\cdots,k_m}^\dag(z)$.
From now on we fix
a lift of Frobenius by $\tilde F_p(t)=t^p$.
By $\mathcal{V}_\mathrm{ KZ}=
\mathcal{O}_X\langle\langle A,B\rangle\rangle$,
we have
$\mathcal{V}^\dag_\mathrm{ KZ}=j^\dag\mathcal{O}_{]\overline{X}_0[}
\langle\langle A,B\rangle\rangle$.

\begin{defn}\label{Bosnia}
The {\sf overconvergent $p$-adic MPL}
$Li_{k_1,\cdots,k_m}^\dag(z)\in 
j^\dag\mathcal{O}_{]\overline{X}_0[}
(]\overline X_0[)$
is the coefficient of $A^{k_m-1}B\cdots A^{k_1-1}B$ in 
$G^\dag_0(z)$ multiplied by $(-1)^m$:
$$
G^\dag_0(z)=1+\cdots+(-1)^mLi_{k_1,\cdots,k_m}^\dag(z)
A^{k_m-1}B\cdots A^{k_1-1}B+\cdots
$$
where $G^\dag_0(z)\in\mathcal{V}^\dag_\mathrm{ KZ}(]\overline X_0[)$
is the image of 
$1\in\tilde F_p^*\mathcal{V}^\dag_\mathrm{ KZ}(]\overline X_0[)=
\mathcal{V}^\dag_\mathrm{ KZ}(]\overline X_0[)=
j^\dag\mathcal{O}_{]\overline{X}_0[}(]\overline X_0[)
\langle\langle A,B\rangle\rangle$
by  \eqref{orange}.
\end{defn}
We note that
since $G^\dag_0(z)$ is group-like we can 
recover and express general coefficients of $G^\dag_0(z)$ 
in terms of $Li_{k_1,\cdots,k_m}^\dag(z)$
as Note \ref{Scottland}.
Our choice of $\tilde F_p$ is a good lifting of the Frobenius
of $\bold P^1_{\bold F_p}\backslash\{1\}$.
So the overconvergent $p$-adic MPL
is analytic on $]0[$ and $]\infty[$, i.e. 
$Li_{k_1,\cdots,k_m}^\dag(z)\in 
j'^\dag\mathcal{O}_{]\overline{X}_0[}(]\overline{X}_0[)$
with
$j':\bold P^1_{\bold F_p}\backslash\{1\}\hookrightarrow{\overline X}_0$
while our $p$-adic MPL $Li_{k_1,\cdots,k_m}(z)$ is a Coleman function,
i.e. belong to the ring $A_\mathrm{ Col}$ of Coleman functions of $X$ 
which contains 
$j^\dag\mathcal{O}_{]\overline{X}_0[}(]\overline{X}_0[)$.
In Lemma \ref{Iceland} we see that overconvergent $p$-adic MPL's
have a different tannakian origin from ours,
so that they are different from ours.
But we have a close relationship between them in 
$A_\mathrm{ Col}\langle\langle A,B\rangle\rangle$:

\begin{thm}\label{Czech}
\begin{equation}\label{pink}
G_0^\dag(A,B)(z)=
G_0(A,B)(z)\cdot
G_0\left(\frac{A}{p} \ ,
\Phi^p_{\mathrm{De}}(A,B)^{-1}\frac{B}{p}\Phi^p_{\mathrm{De}}(A,B)\right)(z^p)^{-1}.  
\end{equation}
\end{thm}

\begin{pf}
Let $z$ be a $\bold Q_p$-valued point of $X$  whose modulo $p$
reduction $z_0$ lies on $X(\bold F_p)$
(If $p=2$, we need to enlarge the base field $\bold Q_p$).
We have a de Rham path $d_z$ in 
$\pi_1^{\mathrm{DR}}(X:\overrightarrow{01},z)(\bold Q_p)$
while we also have a de Rham path $\phi_p(d_{z^p})$
by
$\pi_1^{\mathrm{DR}}(X:\overrightarrow{01},z)\times\bold Q_p
\cong
\pi_1^{p,\mathrm{rig}}(X_0:\overrightarrow{01},z_0)
\cong
\pi_1^{\mathrm{DR}}(X:\overrightarrow{01},z^p)\times\bold Q_p$
in Lemma \ref{Austria}.
Identifying them 
we get a de Rham loop $d_{z}^{-1}\phi_p(d_{z^p})$ in 
$\pi_1^{\mathrm{DR}}(X:\overrightarrow{01})(\bold Q_p)$.

\begin{lem}\label{Iceland}
$i(d^{-1}_{z}\phi_p(d_{z^p}))=G^\dag_0(z)$.
\end{lem}

\begin{pf}
By the categorical equivalence \eqref{red},
we identify 
$\mathcal P_{\overrightarrow{01}}^{\mathrm{DR}}\times \bold Q_p$,
$\mathcal A_{\overrightarrow{01}}^{\mathrm{DR}}\otimes\bold Q_p$ and
$\mathcal V_{\mathrm{KZ}}\otimes\bold Q_p$ with
$\mathcal P_{\overrightarrow{01}}^\dag$,
$\mathcal A_{\overrightarrow{01}}^\dag$ and
$\mathcal V_{\mathrm{KZ}}^\dag$ respectively.
Because $\mathcal P_{\overrightarrow{01}}^{\mathrm{DR}}
=Spec \mathcal A_{\overrightarrow{01}}^{\mathrm{DR}}$,
we have a natural injection
$\mathcal P_{\overrightarrow{01}}^{\mathrm{DR}}\to
(\mathcal A_{\overrightarrow{01}}^{\mathrm{DR}})^\lor$
of pro-objects.
By Lemma \ref{Germany} we get a morphism
$\mathcal P_{\overrightarrow{01}}^{\mathrm{DR}}\to
\mathcal V_{\mathrm{KZ}}$,
whose associated morphism 
$\mathcal P_{\overrightarrow{01}}^\dag\to
\mathcal V_{\mathrm{KZ}}^\dag$
commutes with  Frobenius actions, i.e.
$$
\begin{CD}
\tilde F_p^*\mathcal P_{\overrightarrow{01}}^\dag @>\phi=F_p^{-1}>>
\mathcal P_{\overrightarrow{01}}^\dag \\
@VVV @VVV \\
\tilde F_p^*\mathcal V_{\mathrm{KZ}}^\dag @>\phi=F_p^{-1}>>
\mathcal V_{\mathrm{KZ}}^\dag \\
\end{CD}
$$
commutes.
By taking the fiber $\omega_z$,
we easily see that 
$d_{z^p}\in\pi_1^{\mathrm{DR}}(X:\overrightarrow{01},z^p)(\bold Q_p)
=\omega_{z^p}(\mathcal P_{\overrightarrow{01}}^\dag)
=\omega_z(\tilde F_p^*\mathcal P_{\overrightarrow{01}}^\dag)$
corresponds to 
$1_{(z^p)}\in\omega_{z^p}(\mathcal V^\dag_{\mathrm{KZ}})
=\omega_z(\tilde F_p^*\mathcal V_{\mathrm{KZ}}^\dag)$.
Hence the image 
$\phi_p(d_{z^p})\in
\pi_1^{\mathrm{DR}}(X:\overrightarrow{01},z)(\bold Q_p)
=\omega_z(\mathcal P_{\overrightarrow{01}}^\dag)$
is mapped to 
$G_0^\dag(z)\in\mathcal V_{\mathrm{KZ},(z)}^\dag
=\omega_{z}(\mathcal V_{\mathrm{KZ}}\otimes\bold Q_p)$.
The second row induces morphism  
$\omega_{\overrightarrow{01}}
(\mathcal P_{\overrightarrow{01}}^{\mathrm{DR}})
=\pi_1^{\mathrm{DR}}(X:\overrightarrow{01})(\bold Q_p)\to
\omega_{\overrightarrow{01}}(\mathcal V_{\mathrm{KZ}})\cong
\bold Q_p\langle\langle A,B\rangle\rangle$,
so that the loop
$d_{z}^{-1}\phi_p(d_{z^p})$ corresponds to $G_0^\dag(z)$.
\qed
\end{pf}

By $d^{-1}_{z}\phi_p(d_{z^p})=
d^{-1}_{z}c_{z_0}\cdot \phi_p(d^{-1}_{z^p}c_{z_0})^{-1}$
we get the equality \eqref{pink}
whenever we fix each $\bold Q_p$-point $z$.
By letting $z$ varies for all $\bold Q_p$-points on $]z_0[$
and restricting $G^\dag_0$ and $G_0$ into $]z_0[$,
we see that the equality \eqref{pink} holds in
$\bold Q_p[[t_{z_0}]]\langle\langle A,B \rangle\rangle$
($t_{z_0}$: a local parameter of $]z_0[ \ \bigcap X(\bold Q_p)$),
hence for $A_{Col}|_{]z_0[}$.
By the uniqueness principle of Coleman functions \cite{Col},
the equality \eqref{pink} holds for the whole space.
After all we get Theorem \ref{Czech}.
\qed
\end{pf}

By Lemma \ref{Iceland} we may say that $p$-adic MZV \`{a} la Deligne
is \lq a special value ' of the overconvergent $p$-adic MPL.
By Theorem \ref{Czech} the overconvergent $p$-adic MPL is 
is described very explicitly as a combination of our $p$-adic MPL's.
The following  are the easiest examples.

\begin{eg}\label{Spain}
\begin{enumerate}
\item $Li^\dag_k(z)=Li_k(z)-\frac{1}{p^k}Li_k(z^p)
=\sum\limits_{(n,p)=1}\frac{z^n}{n^k}$.
\item $Li^\dag_{a,b}(z)=
Li_{a,b}(z)-\frac{1}{p^{a+b}}Li_{a,b}(z^p)-
(\frac{1}{p^b}-\frac{1}{p^{a+b}})\zeta_p(a)Li_b(z^p)$
\par
$
-\sum\limits_{r=0}^{a-1}(-1)^r\frac{1}{p^{a-r}}
\binom{b-1+r}{r} Li_{a-r}(z^p)
\{Li_{b+r}(z)-\frac{1}{p^{b+r}}Li_{b+r}(z^p)\}$
\par
$
-(-1)^a\sum\limits_{s=0}^{b-1}\binom{a-1+s}{a-1}
(\frac{1}{p^{b-s}}-\frac{1}{p^{a+b}})
\zeta_p(a+s)Li_{b-s}(z^p)$.
\end{enumerate}
\end{eg}

We note that the coefficient of $A$ of $G_0(z)$ is $\log^az$
while in the overconvergent side the coefficient
of $G_0^\dag(z)$ is $\log^az-\frac{1}{p}\log^az^p=0$
and $Li^\dag_k(z)$ has been studied by Coleman \cite{Col}.
Unver nearly got the same formula (2) in \cite{U} 5.16.
In Example \ref{book}, we will give a Hodge analogue of the
above two formulae. 
Our $p$-adic MPL has log poles around $0$ and $\infty$
however the overconvergent $p$-adic MPL does not.
The formula \eqref{pink} is an algorithm to erase 
log poles at $\infty$ of our $p$-adic MPL.
We remark that the overconvergent $p$-adic MPL is 
rigid analytic on an open rigid analytic subspace $]\overline{X}_0[$
containing $]\bold P^1\backslash\{1\}[$.
 
\begin{rem}\label{Madrid}
The $p$-adic KZ equation
$$
dg=(\frac{A}{z}+\frac{B}{z-1})gdz
$$
does not have a solution with overconvergent function coefficients
but have a solution with Coleman function coefficients 
(for instance $G_0(z)$).
However if we modify it as follows
\begin{equation}\label{Princeton}
dg=\left(\frac{A}{z}+\frac{B}{z-1}\right)gdz
-g\left(\frac{dz^p}{z^p}\frac{A}{p}+
\frac{dz^p}{z^p-1}
\Phi^p_{\mathrm{De}}(A,B)^{-1}\frac{B}{p}\Phi^p_{\mathrm{De}}(A,B)
\right),
\end{equation}
it has a solution with overconvergent function coefficients
(for instance $G^\dag_0(z)$).
We note that
the similar differential equation for the loop
$d_{z^p}^{-1}F_p(d_z)$ is introduced in \cite{U,Y} before.
\end{rem}

\subsection{Hodge side --- analogous results}\label{Hodge}
In this side, we use a Hodge-type comparison
isomorphism (Lemma \ref{Armenia}) between
a de Rham fundamental torsor and a Betti fundamental torsor for $X$.
A tannakian interpretation of MPL's and MZV's are given 
in Proposition \ref{Azerbaijan}
and Proposition \ref{Georgia}.
By using the infinity Frobenius action on a Betti path space
we introduce a new version of MPL's \eqref{Cambridge}.
Its single-valuedness is shown in Theorem \ref{Russia}.
A formula to express them in terms of usual (multi-valued)
MPL's is shown in Proposition \ref{Moldova}.

\begin{nota}\label{Belarus}
Let $z$ be a point in $X(\bold C)$.
We fix a topological path
$b_z\in\pi_1^{\mathrm{top}}(X(\bold C):\overrightarrow{01},z)$.
By abuse of notations we also denote $b_z$ to be the corresponding
Betti path in $\pi_1(X(\bold C):\overrightarrow{01},z)(\bold Q)$.
For $z=\overrightarrow{10}$
we consider a special path 
$b\in\pi_1(X(\bold C):\overrightarrow{01},\overrightarrow{10})(\bold Q)$
which comes from a one-point set 
$\pi_1^{\mathrm{top}}(X(\bold R):
\overrightarrow{01},\overrightarrow{10})$.
\end{nota}

The following Hodge type comparison isomorphism
for fundamental torsors is our basics in this subsection.

\begin{lem}\label{Armenia}
Let $z$ be a $\bold C$-valued point of $X$.
Then we have
\begin{equation}\label{aqua}
\pi_1^{\mathrm{DR}}(X:\overrightarrow{01},z)\times\bold C
\cong
\pi_1^{\mathrm{Be}}(X:\overrightarrow{01},z)\times\bold C.
\end{equation}
\end{lem}

This follows from an equivalence of tannakian categories
\begin{equation}\label{silver}
\mathcal{NC}^\mathrm{DR}(X_{\bold C})\cong
\mathcal{NL}^\mathrm{Be}(X(\bold C))
\end{equation}
(see also \cite{De}).
By identifying the above two fundamental torsors
by \eqref{aqua} 
we get a de Rham loop $d_z^{-1}b_{z}$,
which lies in 
$\pi_1^\mathrm{DR}(X:\overline{01})(\bold C)$.

\begin{prop}\label{Azerbaijan}
Let $z$ be a $\bold C$-valued point of $X$.
The de Rham loop $d_z^{-1}b_{z}$ 
corresponds to the special value at $z$
of the analytic continuation along $b_z$
of the fundamental solution $G_0$
\footnote{
It means the fundamental solution of the KZ equation satisfying
$G_0(z)\approx z^A$ 
constructed by Drinfel'd \cite{Dr}.
}
of the KZ equation by the embedding $i$:
$$
i(d_{z}^{-1}b_{z})=G_0(z).
$$
\end{prop}

\begin{pf}
The proof is given in a similar way to the proof of 
Theorem \ref{Sweden}. 
By \cite{De} \S 12.15-12.16 we have a tannakian interpretation of
horizontal sections of
${\mathcal P}^{\mathrm{DR}}_{\overrightarrow{01}}$,
By Lemma \ref{Germany} we have a tannakian interpretation of 
horizontal sections of
$(\mathcal{V}_\mathrm{ KZ},\nabla_\mathrm{ KZ})$.
Because both the constant term of $G_0(t)$ at $t=0$
and the constant term of $t_0^A$ at $t_0=0$ are equal to $1$,
we have 
$$
b_{\overrightarrow{01},z}(t_0^A|_{\overrightarrow{01}})=
\text{the analytic continuation of $G_0(t)$
along } b_z
$$
which implies
\begin{equation}\label{olive}
d_{z,\Gamma}\cdot b_{\overrightarrow{01},z}\cdot 
d_{\Gamma,\overrightarrow{01}}(1)=G_0(z).
\end{equation}
By \eqref{brown}, \eqref{olive} and
$d_{z,\Gamma}=d_{\overrightarrow{01},\Gamma}\cdot d_{z}^{-1}$,
we get Proposition \ref{Azerbaijan}.
\qed
\end{pf}

In each coefficient of $G_0(z)$
there appears MPL \eqref{purple}.
In precise the coefficient of $A^{k_m-1}B\cdots A^{k_1-1}B$ in $G_0(z)$
is $(-1)^mLi_{k_1,\cdots,k_m}(z)$:
$$
G_0(z)=1+\log zA-Li_1(z)B+\cdots+
(-1)^mLi_{k_1,\cdots,k_m}(z)A^{k_m-1}B\cdots A^{k_1-1}B+\cdots.
$$
By Proposition \ref{Azerbaijan} we may say that the loop $d_z^{-1}b_{z}$
is a tannakian origin of MPL's.
We remark that a choice of the Betti path $b_z$
corresponds to a choice of each
branch of analytic continuation of MPL. 

For $z=\overrightarrow{10}$ similarly we have

\begin{prop}\label{Georgia}
The de Rham loop 
$d^{-1}b\in\pi_1^\mathrm{DR}(X:\overrightarrow{01})(\bold C)$
corresponds to the Drinfel'd associator 
$\Phi_\mathrm{KZ}\in\bold C\langle\langle A,B\rangle\rangle$
by the embedding $i$, i.e.
$$
i(d^{-1}b)=\Phi_\mathrm{KZ}.
$$
\end{prop}

\begin{pf}
The Drinfel'd associator $\Phi_\mathrm{KZ}$
is the series constructed in \cite{De}, 
which is also equal to the limit value 
$$
{\underset{\epsilon\to 0}{\lim}} 
\exp\left(-\log\epsilon\cdot B\right)\cdot 
G_0(1-\epsilon)
\in\bold C\langle\langle A,B\rangle\rangle
$$
as \eqref{rose}.
In other word it is 
the constant term of $G_0(z)$ at $z=1$.
Then by a similar argument to the proof of Theorem \ref{Finland}
we have 
$b_{\overrightarrow{01},\overrightarrow{10}}
(t_0^A|_{\overrightarrow{01}})
=\Phi_\mathrm{KZ}\cdot
t_1^B|_{\overrightarrow{10}}$,
whence
\begin{equation}\label{white}
d_{\overrightarrow{10},\Gamma}\cdot b\cdot d_{\Gamma,\overrightarrow{01}}(1)
=\Phi_\mathrm{KZ}.
\end{equation}
By \eqref{brown}, \eqref{white} and
$d_{\overrightarrow{10},\Gamma}=d_{\overrightarrow{01},\Gamma}\cdot d^{-1}$,
we prove Proposition \ref{Georgia}.
\qed
\end{pf}

In each coefficient of 
the Drinfel'd associator $\Phi_\mathrm{KZ}$ 
there appears MZV \eqref{gold}.
In precise the coefficient of $A^{k_m-1}B\cdots A^{k_1-1}B$ in 
$\Phi_\mathrm{KZ}$ 
is $(-1)^m\zeta(k_1,\cdots,k_m)$:
$$
\Phi_\mathrm{KZ}=1+\cdots+
(-1)^m\zeta(k_1,\cdots,k_m)A^{k_m-1}B\cdots A^{k_1-1}B+\cdots.
$$
By Proposition \ref{Georgia} we may say that the loop $d^{-1}b$
is a tannakian origin of MZV's.

The above proposition enables us to calculate the period map
\begin{equation}\label{enable}
p:\pi_1^\mathrm{ Be}(X(\bold C):\overrightarrow{01})(\bold C)\to
\pi_1^\mathrm{ DR}(X:\overrightarrow{01})(\bold C).
\end{equation}

In order to do that,
as in \S \ref{Berthelot-Ogus} we fix generator of
$\pi_1^\mathrm{ Be}(X(\bold C):\overrightarrow{01})(\bold Q)$
and its parameterization.

\begin{nota}\label{Chile}
We denote the loop in 
$\pi_1^\mathrm{ Be}(X(\bold C):\overrightarrow{01})(\bold Q)$
(resp. in 
$\pi_1^\mathrm{ Be}(X(\bold C):\overrightarrow{10})(\bold Q)$)
which comes from a generator in 
$\pi_1^\mathrm{ top}(T_0^\times(\bold C):\overrightarrow{01})$
(resp. in 
$\pi_1^\mathrm{ top}(T_1^\times(\bold C):\overrightarrow{10})$)
going around the point $0\in T_0(\bold C)$ (resp. $0\in T_1(\bold C)$)
counterclockwisely by $x_\mathrm{ Be}$ (resp. $y'_\mathrm{ Be}$).
Put $y_\mathrm{ Be}=b^{-1}y'_\mathrm{ Be}b$
(for $b$ see Notation \ref{Belarus}).
Easily we see that two loops $x_\mathrm{ Be}$ and $y_\mathrm{ Be}$
are free generators of the pro-algebraic group 
$\pi_1^\mathrm{ Be}(X(\bold C):\overrightarrow{01})(\bold Q)$.
Wojtkowiak \cite{W} consider the embedding
\begin{equation}\label{tan}
j:\pi_1^\mathrm{ Be}(X(\bold C):\overrightarrow{01})(\bold Q)
\hookrightarrow\bold Q\langle\langle A,B\rangle\rangle
\end{equation}
which sends 
$x_\mathrm{ Be}\mapsto e^A:=1+\frac{A}{1!}+\frac{A^2}{2!}+\cdots$ and
$y_\mathrm{ Be}\mapsto e^B:=1+\frac{B}{1!}+\frac{B^2}{2!}+\cdots$.
We note that this morphism is isomorphic to the embedding of
$\pi_1^\mathrm{Be}(X(\bold C):\overrightarrow{01})(\bold Q)$
into the universal embedding algebra of 
$Lie\pi_1^\mathrm{Be}(X(\bold C):\overrightarrow{01})(\bold Q)$.
\end{nota}

We note that this embedding $j$ is not compatible with $i$ 
in \eqref{navy} under the Hodge comparison isomorphism \eqref{aqua}.
As is similar to the proof of Theorem \ref{Netherland},
we extend the period map \eqref{enable} to
$$p:\bold C\langle\langle A,B\rangle\rangle\to 
\bold C\langle\langle A,B\rangle\rangle$$
via the embedding $i$ \eqref{navy} and $j$ \eqref{tan}.

\begin{lem}\label{Vancouver}
$p(A)=2\pi iA$, $p(B)=\Phi_{\mathrm{KZ}}^{-1}(2\pi iB)\Phi_{\mathrm{KZ}}$.
\end{lem}

\begin{pf}
By the compatibility of the period map on the tangential morphisms
$\bold Q_\mathrm{ Be}(1)\overset{\exp}{=}
\pi_1^\mathrm{Be}(T^\times_0(\bold C):\overrightarrow{01})
\to
\pi_1^\mathrm{Be}(X(\bold C):\overrightarrow{01})$
and
$\bold Q_\mathrm{ DR}(1)\overset{\exp}{=}
\pi_1^\mathrm{DR}(T^\times_0:\overrightarrow{01})
\to
\pi_1^\mathrm{DR}(X:\overrightarrow{01})$,
we have $p(x_\mathrm{ Be})=x^{2\pi i}$.
Similarly we have $p(y'_\mathrm{ Be})=y'^{2\pi i}$.
Since the period map is compatible with the torsor structure
we have
$p(y_\mathrm{ Be})=p(b^{-1}y'_\mathrm{ Be}b)=b^{-1}y'^{2\pi i}b
=(d^{-1}b)^{-1}y^{2\pi i}d^{-1}b$.
By Proposition \ref{Georgia}, we get the claim.
\qed
\end{pf}

Our next work is on a description of a Frobenius action.
To begin with we consider the infinity Frobenius action $\phi_\infty$ on 
$\pi_1^{\mathrm{DR}}(X:\overrightarrow{01},\overrightarrow{10})
(\bold C)$
by transmitting the infinity Frobenius action on 
$\pi_1^{\mathrm{Be}}(X(\bold C):\overrightarrow{01},\overrightarrow{10})
(\bold C)$ by Lemma \ref{Armenia}.
It gives a new path $\phi_\infty(d)$ in 
$\pi_1^{\mathrm{DR}}(X:\overrightarrow{01},\overrightarrow{10})
(\bold C)$.
We denote $\Phi^-_{\mathrm{KZ}}$ to be the series in
$\bold C\langle\langle A,B\rangle\rangle$
which corresponds to
$d^{-1}\phi_\infty(d)$:
$$
i(d^{-1}\phi_\infty(d))=\Phi^-_{\mathrm{KZ}}.
$$
It is a Hodge counterpart of the $p$-adic Deligne associator
$\Phi^p_{\mathrm{De}}$.
By same arguments to the proof of 
Lemma \ref{Ireland} and Theorem \ref{Netherland}
the extension of the infinity Frobenius action on 
$\pi_1^{\mathrm{DR}}(X:\overrightarrow{01})(\bold C)$
into $\bold C\langle\langle A,B\rangle\rangle$
by $i$ is described as follows:

\begin{lem}\label{Kylgyz}
$\phi_\infty(A)=-A$, \ $\phi_\infty(B)=
{\Phi^-_{\mathrm{KZ}}}^{-1}(-B)\Phi^-_{\mathrm{KZ}}$.
\end{lem}

This gives the following formula analogous to \eqref{Amsterdam}.

\begin{lem}\label{Kazakhstan}
$$
\Phi_{\mathrm{KZ}}(A,B)= \Phi^-_{\mathrm{KZ}}(A,B)\cdot
\Phi_{\mathrm{KZ}}\left(-A,
\Phi^-_{\mathrm{KZ}}(A,B)^{-1}(-B)\Phi^-_{\mathrm{KZ}}(A,B)\right).  
$$
\end{lem}

This formula looks nasty.
But according to the convention in \cite{F0}
of the Grothendieck-Teichm\"{u}ller group \cite{Dr}
we rewrite simpler
$$
(-1,\Phi^-_{\mathrm{KZ}})=(-1,\Phi_{\mathrm{KZ}})\circ
(1,\Phi_{\mathrm{KZ}})^{-1}. 
$$
Let $z$ be a $\bold C$-valued point of $X$.
By a same argument to above,
we get a de Rham loop
$d_{z}^{-1}\phi_\infty(d_{\bar z})\in
\pi_1^\mathrm{ DR}(X:\overline{01})(\bold C)$.
Put
$G_0^-(z):=i(d_{z}^{-1}\phi_\infty(d_{\bar z}))$.
The following is an analogue to Theorem \ref{Czech}.

\begin{prop}\label{Moldova}
\begin{equation}\label{gray}
G_0^-(A,B)(z)=
G_0(A,B)(z)\cdot
G_0\left(-A,\Phi^-_{\mathrm{KZ}}(A,B)^{-1}(-B)
\Phi^-_{\mathrm{KZ}}(A,B)\right)(\bar z)^{-1}.  
\end{equation}
\end{prop}

\begin{pf}
By $d^{-1}_{z}\phi_\infty(d_{\bar z})=
d^{-1}_{z}b_{z}\cdot \phi_\infty(d^{-1}_{\bar z}\overline{b_{z}})^{-1}$
and Lemma \ref{Kylgyz} we get the equality.
\qed
\end{pf}

We introduce a new version of MPL $Li^-_{k_1,\cdots,k_m}(z)$
to be the coefficient of 
$A^{k_m-1}B$ $\cdots$ $A^{k_1-1}B$ in $G^-_0(z)$
multiplied by $(-1)^m$:
\begin{equation}\label{Cambridge}
G^-(z)=:1+\cdots+(-1)^mLi^-_{k_1,\cdots,k_m}(z)
A^{k_m-1}B\cdots A^{k_1-1}B+\cdots.
\end{equation}

We note that since $G^-(z)$ is group-like we can 
recover and express general coefficients of $G^-(z)$ 
in terms of $Li^-_{k_1,\cdots,k_m}(z)$ and $\log |z|$
as Note \ref{Scottland}.

\begin{thm}\label{Russia}
The coefficient $Li^-_{k_1,\cdots,k_m}(z)$ is single-valued 
and real-analytic
on $X(\bold C)=\bold P^1(\bold C)\backslash\{0,1,\infty\}$.
\end{thm}

\begin{pf}
The real-analyticity follows from \eqref{gray}
because both $G_0(z)$ and $G_0(\bar z)$ are real-analytic.
The single-valuedness is immediate by definition.
\qed
\end{pf}

The single-valued MPL's are not equal to usual MPL's
because they have a different tannakian origins.
But by Proposition \ref{Moldova} the single-valued MPL  
is described very explicitly as a combination of usual MPL's.
The following  are the easiest examples:

\begin{eg}\label{book}
\begin{enumerate}
\item $Li^-_k(z)=Li_k(z)-\sum\limits_{a=0}^{k-1}(-1)^{k-a}
\frac{(\log|z|^2)^a}{a!}Li_{k-a}(\bar z)$
\item $Li^-_{a,b}(z)=
Li_{a,b}(z)-\sum\limits_{r=0}^{a-1}\sum\limits_{s=0}^r
[(-1)^{a+r+s}\binom{b-1+s}{s}
\frac{(\log |z|^2)^{r-s}}{(r-s)!}Li_{a-r}(\bar z)\cdot$
\par
$
\{
Li_{b+s}(z)-\sum\limits_{w=0}^{b+s-1}
(-1)^{b+s+w}\frac{(\log |z|^2)^w}{w!}Li_{b+s-w}(\bar z)
\}]
-\sum\limits_{u=0}^{b-1}\frac{(\log |z|^2)^u}{u!}\cdot
$
\par
$
[(-1)^{a+b+u}Li_{a,b-u}(\bar z)+
\{(-1)^{b+u}-(-1)^{a+b+u}\}\zeta(a)Li_{b-u}(\bar z)
$
\par
$
+\sum\limits_{v=0}^{b-u-1}\{(-1)^{a+b+u+v}-(-1)^{b+u}\}
\binom{a+v-1}{a-1}\zeta(a+v)Li_{b-u-v}(\bar z)]$.
\end{enumerate}
\end{eg}

The coefficient of $A$ of $G_0(z)$ is $\log z$
while the coefficient of $G_0^-(z)$ is $\log\bar z+\log z=\log |z|^2$.
We also remark that Zhao \cite{Zh} constructed 
a single-valued version of several variable MPL
$Li_{k_1,\cdots,k_m}(z_1,\cdots,z_m)=\sum\limits_{0<n_1<\cdots<n_m}
\frac{z_1^{n_1}\cdots z_m^{n_m}}{n_1^{k_1}\cdots n_m^{k_m}}$.

\begin{rem}\label{Ukraine}
In \cite{Za} Zagier studied another variant of polylogarithm
$$
P_k(z)=\Re_k\sum\limits_{a=0}^{k-1}\frac{B_a}{a!}(\log|z|^2)^a
Li_{k-a}(z)
$$
where $\Re_k$ denote $\Re$ or $\Im$ depending whether
$k$ is odd or even and
$\frac{t}{e^t-1}=\sum\limits_{n=0}^\infty B_n\frac{t^n}{n!}$. 
He showed that its single-valuedness and real-analyticity.
Beilinson and Deligne \cite{BD} give its interpretations in terms of 
variants of mixed Hodge structures.
Their computations \cite{BD} \S 1.5 says that
$P_k(z)$ appears as the coefficient of $A^{k-1}B$ of 
$\log G^-_0(z)$ multiplied by $-\frac{1}{2}$,
whence the tannakian origin of $P_k(z)$ is also 
$d^{-1}_{z}\phi_\infty(d_{\bar z})$.
We notice that $P_k(z)$ and $Li_k^-(z)$ are related with each other by
\begin{equation}\label{Oxford}
P_k(z)=\frac{1}{2}\sum\limits_{i=0}^{k-1}\frac{B_i}{i!}\{\log |z|^2\}^i
Li_{k-i}^-(z)
\end{equation}
and
\begin{equation*}
Li_k^-(z)=2\sum\limits_{i=0}^{k-1}\frac{(\log |z|^2)^i}{(i+1)!}P_{k-i}(z).
\end{equation*}
\end{rem}

Our MPL has monodromies around $0$, $1$ and $\infty$
however the variant does not.
The formula \eqref{gray} gives an algorithm to erase 
monodromies of our MPL at $0$, $1$ and $\infty$.
 
\begin{rem}
The KZ equation
$$
dg=(\frac{A}{z}+\frac{B}{z-1})gdz
$$
have a solution with multi-valued complex analytic function coefficients 
(for instance $G_0(z)$) 
whereas the following modification analogous to \eqref{Princeton}
$$
dg=\left(\frac{A}{z}+\frac{B}{z-1}\right)gdz
-g\left(\frac{d\bar z}{\bar z}(-A)+
\frac{d\bar z}{\bar z-1}
\Phi^-_{\mathrm{KZ}}(A,B)^{-1}(-B)\Phi^-_{\mathrm{KZ}}(A,B)
\right)
$$
does have a solution with single-valued real-analytic 
function coefficients(for instance $G^-_0(z)$).
\end{rem}

\subsection{Artin side --- reviews}\label{Artin}
In this side we will discuss $l$-adic \'{e}tale  analogues of
MPL's and MZV's which we call $l$-adic MPL's and 
$l$-adic multiple Soul\'{e} elements respectively and
whose definitions are imitations of 
Wojtkowiak's $l$-adic polylogarithms \cite{W}.\par

Let $l$ be a prime (we do not assume whether $l$ is equal to the
previous prime $p$ in \S \ref{Berthelot-Ogus} or not).
Let $\overline{\bold Q}$ be 
the algebraic closure of $\bold Q$ in $\bold C$.
The following Artin-type comparison isomorphism take place of 
Lemma \ref{Austria} and Lemma \ref{Armenia}.

\begin{lem}\label{Argentina}
Suppose that $x$ and $y$ are $\overline{\bold Q}$-valued points 
of $X$ or $T^\times_s$ ($s\in\{0,1,\infty\}$)
then there exists an isomorphism 
$$
\pi_1^\mathrm{ Be}(X(\bold C):x,y)\times_\bold Q\bold Q_l\cong
\pi_1^{l,\text{\'{e}t}}(X_{\overline{\bold Q}}:x,y).
$$
\end{lem}

It follows from the comparison isomorphism \cite{SGA1} XII Corollary 5.2
(see also \cite{De} \S 13.11).

Let
$i:\pi_1^{l,\text{\'{e}t}}(X_{\overline{\bold Q_l}}:\overrightarrow{01})
(\bold Q_l)\hookrightarrow\bold Q_l\langle\langle A,B\rangle\rangle$
be the embedding associated with $j$ \eqref{tan}.
Suppose that $z\in X(\bold Q)$ and
$b_z\in\pi_1^\mathrm{ top}(X(\bold C):\overrightarrow{01},z)$.
It determines a Betti path in 
$\pi_1^\mathrm{ Be}(X(\bold C):\overrightarrow{01},z)(\bold Q)$
and we regard this to be an $l$-adic \'{e}tale path
(which we denote $b_z$ by abuse of notation) in 
$\pi_1^{l,\text{\'{e}t}}(X_{\overline{\bold Q}}:\overrightarrow{01},z)
(\bold Q_l)$ Lemma \ref{Argentina}.
Suppose that $\sigma\in Gal({\overline{\bold Q}}/{\bold Q})$.
We obtain another path $\sigma(b_z)$ by the Galois group action on 
$\pi_1^{l,\text{\'{e}t}}(X_{\overline {\bold Q}}:\overrightarrow{01},z)(\bold Q_l)$. 
Combining $\sigma(b_z)$ with the inverse $b_z^{-1}$,
we get an $l$-adic \'{e}tale loop
$b_z^{-1}\sigma(b_z)\in
\pi_1^{l,\text{\'{e}t}}(X_{\overline {\bold Q}}:\overrightarrow{01})
(\bold Q_l)$. 
Wojtkowiak \cite{W} introduced 
$l$-adic polylogarithm $\ell_k^z(\sigma)$ ($k\geqslant 1$) 
to be the coefficient of $A^{k-1}B$ of 
$\log j(b_z^{-1}\sigma(b_z))$ multiplied by $(-1)^{k-1}$ (cf. \cite{NW}).
In accordance with expressions in \S \ref{Berthelot-Ogus} and
\S \ref{Hodge},
we consider the following $l$-adic \'{e}tale analogues of
MPL's and MZV's.

\begin{defn}\label{Colombia}
Let $m,k_1,\cdots,k_{m}\geqslant 1$, $z\in X(\bold Q)$ and 
$\sigma\in Gal(\overline{\bold Q}/\bold Q)$.
The {\sf $l$-adic multiple polylogarithm} 
(shortly the {\sf $l$-adic MPL}) $Li^l_{k_1,\cdots,k_m}(z)(\sigma)$
denotes the coefficient of $A^{k_m-1}B\cdots A^{k_1-1}B$ in 
$j(b_z^{-1}\sigma(b_z))$ multiplied by $(-1)^m$
and {\sf $l$-adic multiple Soul\'{e} element}
(shortly the {\sf $l$-adic MSE})
$\zeta_l(k_1,\cdots,k_m)(\sigma)$ ($k_m>1$)
denotes the coefficient of $A^{k_m-1}B\cdots A^{k_1-1}B$ in 
the {\sf $l$-adic Ihara associator}
$\Phi^l_\sigma(A,B)=j(b^{-1}\sigma(b))$ multiplied by $(-1)^m$.
\end{defn}

By making a same computations similar to \cite{BD} \S 1.5,
we have
\begin{equation}\label{darkgray}
\ell_k^z(\sigma)=(-1)^{k}\sum\limits_{i=0}^{k-1}
\frac{B_i}{i!}\{\rho_z(\sigma)\}^i
Li^l_{k-i}(z)(\sigma).
\end{equation}
Here $\rho_z(\sigma)$ is the generalized $l$-adic Kummer 1-cocycle
in \cite{NW} Definition 3.
By comparing it with \eqref{Oxford},
we may say that $Li^l_k(z)(\sigma)$
is a $l$-adic analogue of the usual (Leibnitz's) polylogarithm
$Li_k(z)=\sum\limits_{n=1}^{\infty}\frac{z^n}{n!}$
and Wojtkowiak's $\ell_k^z(\sigma)$ is a $l$-adic analogue of
the Beilinson-Deligne's polylogarithm $D_k(z)$.
By asking for $p$-adic Hodge theory, we may give new $p$-adic MPL's
valued on the ring of $B_\mathrm{ crys}$ of Fontaine's $p$-adic periods,
which might unify $p$-adic MPL's in \S \ref{Berthelot-Ogus}
and those above (for $l=p$).

The following may explain why we call 
the multiple Soul\'{e} elements.

\begin{eg}[\cite{I90} \S 6.3.Theorem]\label{Ecuador}
\begin{enumerate}
\item When $k$ is odd,
$\zeta_l(k)(\sigma)=
\frac{\kappa^{(l)}_k(\sigma)}{(1-l^{k-1})\dot (k-1)!}$.
Here $\kappa^{(l)}_k(\sigma)$ is  $l$-adic Soul\'{e} element
(see \cite{I90} \S 6.2 for its definition).
\item When $k$ is even,
$\zeta_l(k)(\sigma)=\frac{B_k}{2(k!)}\{1-\chi_l(\sigma)^k\}$.
Here $B_k$ is the Bernoulli number defined by 
$\frac{t}{e^t-1}=\sum\limits_{n=0}^{\infty}B_n\frac{t^n}{n!}$
and $\chi_l(\sigma)$ is the $l$-adic cyclotomic character.
\end{enumerate}
\end{eg}

The author expects that MSE is a kind of \lq multiple Euler system'
which helps our understanding of multiple zeta functions.\par

We end this section by giving a brief observation on 
$l$-adic behavior of $\zeta_l(k)(\sigma)$ and $p$-adic behavior of
$\zeta_p(k)$ with respect to weights.
 
\begin{rem}\label{Uruguai}
Let $\sigma\in Gal(\overline{\bold Q}/\bold Q)$.
By the expression 
$\kappa^l_k(\sigma)=\int_{\bold Z_l^\times}x^{k-1}d\mu^\sigma$
($\mu^\sigma$: the measure associated with the Kummer
distribution \cite{NW} \S 2) 
we see that $(k-1)!(1-l^{k-1})\zeta_l(k)(\sigma)$ admits 
a nice $l$-adic behavior with respect to $k$,
that is, for $k\equiv k'\mod (l-1)l^M$ ($M\in\bold N$) we have
\begin{equation}\label{lightgray}
(k-1)!(1-l^{k-1})\zeta_l(k)(\sigma)\equiv 
(k'-1)!(1-l^{k'-1})\zeta_l(k')(\sigma) 
\mod l^{M+1}\bold Z_l.
\end{equation}
On the other hand in the $p$-adic setting, by 
$L_p(1-k,\omega^{1-k})=-(1-p^{k-1})\frac{B_k}{k}=(1-p^{k-1})\zeta(1-k)$
and $L_p(k,\omega^{1-k})=(1-\frac{1}{p^k})\zeta_p(k)$
we see that $(1-\frac{1}{p^k})\zeta_p(k)$ admits a nice
$p$-adic behavior with respect to $k$,
that is, for $k\equiv k'\mod (p-1)p^M$ ($M\in\bold N$) we have
\begin{equation}\label{whitesmoke}
(1-\frac{1}{p^k})\zeta_p(k)\equiv 
(1-\frac{1}{p^{k'}})\zeta_p(k') 
\mod p^{M+1}\bold Z_p.
\end{equation}
\end{rem}

The author is not sure whether we have \lq multiple analogues' of 
\eqref{lightgray} and \eqref{whitesmoke}.
He is grateful to the referee who suggested him that
Coleman's integral formula (\cite{Col} Lemma 7.2)
$Li^\dag_k(z)=\int_{\bold Z^\times_p}x^{-k}d\mu_z(x)$
for $z\in\bold C_p\backslash\{|z-1|_p<1\}$ might help
look for such relations.
Here $\mu_z$ is a measure on $\bold Z_p$ given by
$\mu_z(a+p^n\bold Z_p)=\frac{z^a}{1-z^{p^n}}$ 
with $n\in\bold N$ and $0\leqslant a< p^n$.

\section{Motivic views}\label{motive}
This section is complementary.
We will consider three boxes related to 
the algebra generated by $p$-adic MZV's. 
In \S \ref{Drinfeld} our special box is
Drinfel'd's  \cite{Dr}
pro-algebraic bi-torsor where the Grothendieck-Teichm\"{u}ller
pro-algebraic groups act.
In \S \ref{Racinet} our special box is 
Racinet's \cite{R} pro-algebraic torsor 
defined by double shuffle relations.
In \S \ref{Goncharov} our special box is 
pro-algebraic bi-torsor of 
Deligne-Goncharov's \cite{DG} motivic Galois group.
The relationship between these three box will be discussed.

\subsection{Drinfel'd's context}\label{Drinfeld}
In this subsection we will add the Berthelot-Ogus part
into the story of \cite{F0} owing Unver's result \cite{U1}
and discuss special elements in Ihara's stable derivation algebra
which are $p$-adic analogues of the Drinfel'd's element $\tilde\psi$
in \cite{Dr} Proposition 6.3.

In his celebrated paper \cite{Dr}, Drinfel'd introduced 
the Grothendieck-Teichm\"{u}ller torsor,
a triple $(\underline{GRT}_1,\underline{M}_1,\underline{GT}_1)$,
where $\underline{M}_1$ is the pro-algebraic bi-torsor with
a left (resp. right) action of the Grothendieck-Teichm\"{u}ller
pro-algebraic group $\underline{GRT}_1$ (resp. $\underline{GT}_1$)
(for their presentations see also \cite{F0}\S 2).
In \cite{F0}\S 4 we explained that for a prime $l$
the $l$-adic Galois image pro-algebraic group 
$\underline{Gal}^{(l)}_{\bold Q_l}$
(which is defined in loc.cit.
to be the Zariski closure of the image of 
the Galois group 
$Gal(\overline{\bold Q}/\bold Q(\mu_{l^\infty}))$-action
in $\underline{Aut}\pi_1^{l,{\text{\'{e}t}}}
(X_{\overline{\bold Q}}:\overrightarrow{01})(\bold Q_l)$)
lies in $\underline{GT}_1\times\bold Q_l$, i.e.
we have
\begin{equation}\label{apple}
\varPhi^{(l)}_{\bold Q_l}:
\underline{Gal}^{(l)}_{\bold Q_l}\hookrightarrow
\underline{GT}_1\times\bold Q_l,
\end{equation}
whereas in \cite{F0} \S 3 we also saw that 
the spectrum of the $\bold Q$-algebra 
$\frak Z (\subset\bold C)$ generated by 
$(2\pi i)^{-(k_1+\cdots+k_m)}\zeta(k_1,\cdots,k_m)$
($m,k_i\in\bold N$, $k_m>1$)
is embedded into $\underline{M}_1$, i.e. we have
\begin{equation}\label{peach}
\varPhi_\mathrm{ Hod}:
Spec \frak Z\hookrightarrow \underline{M}_1.
\end{equation}
We note that \eqref{apple} and \eqref{peach} are expected 
to be isomorphism \cite{F0}.
In \cite{F1} we announced to make up still lacking 
Berthelot-Ogus part in our next upcoming paper
\lq $p$-adic multiple zeta values III';
we planned to construct an embedding
\begin{equation}\label{melon}
\varPhi^{(p)}_\mathrm{ crys}:
Spec \frak Z^{(p)}\hookrightarrow \underline{GRT}_1
\end{equation}
where ${\frak Z}^{(p)}$ is a subalgebra of $\bold Q_p$ 
generated by all $p$-adic MZV's.
But we would cancel it and give an explanation below instead
because a recent Unver's paper \cite{U1} nearly achieved it.

\begin{prop}\label{Japan}
The $p$-adic Drinfel'd associator $\Phi^p_\mathrm{ KZ}(A,B)$
satisfies the following defining equations  
of $\underline{GRT}_1$:
$$
\begin{cases}
{\rm(0)} \ \   \Phi^p_\mathrm{ KZ}
\in exp[\Bbb L_{\bold Q_p}^{\land},\Bbb L_{\bold Q_p}^{\land}]\\
{\rm(i)} \ \   \Phi^p_\mathrm{ KZ}(A,B)
\Phi^p_\mathrm{ KZ}(B,A)=1  \\
{\rm(ii)} \ \  \Phi^p_\mathrm{ KZ}(C,A)\Phi^p_\mathrm{ KZ}
(B,C)\Phi^p_\mathrm{ KZ}(A,B)=1
 \ \ \text{for} \ \ A+B+C=0 \\
{\rm(iii)}  \ \ 
\Phi^p_\mathrm{ KZ}(X_{1,2} ,X_{2,3})
\Phi^p_\mathrm{ KZ}(X_{3,4} ,X_{4,5})
\Phi^p_\mathrm{ KZ}(X_{5,1} ,X_{1,2})\cdot \\
\qquad\qquad\qquad
\Phi^p_\mathrm{ KZ}(X_{2,3} ,X_{3,4})
\Phi^p_\mathrm{ KZ}(X_{4,5} ,X_{5,1})=1 \qquad
\text{in} \qquad  \underline{U\frak P^{(5)}_\centerdot}(\bold Q_p).
\end{cases}
$$
Here $\Bbb L_{\bold Q_p}^{\land}$ stands for 
the completed free Lie algebra  generated by $A$ and $B$ and
$\underline{U\frak P^{(5)}_\centerdot}$ stands for 
the universal enveloping algebra of the pure sphere 5-braid Lie algebra
$\frak P^{(5)}$ with standard generators
$X_{i,j}$ ($1\leqslant i,j\leqslant 5$).
\end{prop}

\begin{pf}
By
$(p,\Phi^p_\mathrm{ De})=(p,1)\circ(1,\Phi^p_\mathrm{ De}(pA,pB))$
and 
$(p,\Phi^p_\mathrm{ KZ})=(p,1)\circ(1,\Phi^p_\mathrm{ KZ}(pA,pB))$
in our convention of $\underline{GRT}_1$, 
\eqref{teal} is rewritten into
\begin{equation}\label{pie}
(1,\Phi^p_\mathrm{ De}(pA,pB))=(1,\Phi^p_\mathrm{ KZ}(pA,pB))\circ
(1,\Phi^p_\mathrm{ KZ}(A,B))^{-1}.
\end{equation}
In \cite{U1} Unver showed that the $p$-adic Deligne associator
$\Phi^p_\mathrm{ De}(A,B)$ satisfies the above relations,
in other words
$(1,\Phi^p_\mathrm{ De}(A,B))\in\underline{GRT}_1(\bold Q_p)$.
Thus
$(1,\Phi^p_\mathrm{ De}(p^nA,p^nB))\in\underline{GRT}_1(\bold Q_p)$
for $n\in\bold N$.
Put 
$$
(1,a_n)=(1,\Phi^p_\mathrm{ De}(pA,pB))^{-1}\circ
(1,\Phi^p_\mathrm{ De}(p^2A,p^2B))^{-1}\circ\cdots\circ
(1,\Phi^p_\mathrm{ De}(p^nA,p^nB))^{-1}.
$$
Because $\underline{GRT}_1$ is closed under the multiplication,
$(1,a_n)$ must belongs to $\underline{GRT}_1(\bold Q_p)$.
By \eqref{pie},
$$
(1,\Phi^p_\mathrm{ De}(p^nA,p^nB))=
(1,\Phi^p_\mathrm{ KZ}(p^nA,p^nB))\circ
(1,\Phi^p_\mathrm{ KZ}(p^{n-1}A,p^{n-1}B))^{-1}.
$$
So we have
$(1,a_n)=(1,\Phi^p_\mathrm{ KZ}(A,B))\circ
(1,\Phi^p_\mathrm{ KZ}(p^nA,p^nB))^{-1}$.
As $n$ goes $\infty$,
$\Phi^p_\mathrm{ De}(p^nA,p^nB)$ converges to $1$ in
$\bold Q_p\langle\langle A,B\rangle\rangle$
which is a projective limit of finite dimensional topological 
$\bold Q_p$-vector space.
Hence
the series $a_n$ ($n\geqslant 1$) should converge to 
$\Phi^p_\mathrm{ KZ}(A,B)$
in $\bold Q_p\langle\langle A,B\rangle\rangle$.
By the locally-compactness of the topological group 
$\underline{GRT}_1(\bold Q_p)$ 
the limit
$(1,\Phi^p_\mathrm{ KZ}(A,B))$ must lie on 
$\underline{GRT}_1(\bold Q_p)$,
whence we get the claim.
\qed
\end{pf}

By Proposition \ref{Japan} we get a morphism \eqref{melon}.
The injectivity is easy to check. \par
Since there appear $p$-adic MZV's in each coefficient of 
the $p$-adic Drinfel'd associator $\Phi^p_\mathrm{ KZ}$,
we get algebraic relations among them.
It is remarkable that we get  
$\zeta_p(2k)=0$ ($k\geqslant1$) 
from (0), (i) and (ii)
(proved in the same way to \cite{De}\S 18.16), that is, 
we have a geometric proof for $L_p(2k,\omega^{1-2k})=0$.\par
In \cite{F0} \S 6 by discussing a weight filtration on our torsor,
we associated \eqref{apple} and \eqref{peach} with two surjections
between their coordinate rings,
$Gr\varPhi^{(l)\sharp}_{\bold Q_l}:
{\mathcal O}(\underline{GRT}_1)_\centerdot\widehat\otimes\bold Q_l
\twoheadrightarrow
Gr^W_\centerdot{\mathcal O}(\underline{Gal}^{(l)}_{\bold Q_l})$
and
$Gr\varPhi_\mathrm{ Hod}^\sharp:
{\mathcal O}(\underline{GRT}_1)_\centerdot\twoheadrightarrow
\Bigl(Z_\centerdot\bigl/(\pi^2)Z_\centerdot\Bigr)$
where
$Z_\centerdot$
is a graded $\bold Q$-algebra 
whose degree $w$-component ($w\geqslant 1$) is 
a $\bold Q$-vector space $Z_w(\subset\bold R)$ 
generated by all MZV's with weight $w$
and $Z_0=\bold Q$.
By \eqref{melon}, we also obtain
$$
\varPhi_\mathrm{ crys}^{(p)\sharp}:
{\mathcal O}(\underline{GRT}_1)_\centerdot\twoheadrightarrow 
Z_\centerdot^{(p)}
$$
where $Z^{(p)}_\centerdot$ 
is a graded $\bold Q$-algebra 
whose degree $w$-component is a $\bold Q$-vector space 
$Z^{(p)}_w(\subset\bold Q_p)$ 
generated by all $p$-adic MZV's with weight $w$ and
$Z^{(p)}_0=\bold Q$.

\begin{rem}\label{Korea}
\begin{enumerate}
\item
In \cite{F0}\S 6.6.1, we introduced an invertible element
$\Phi_{GRT}$ in 
${\mathcal O}(\underline{GRT}_1)\langle\langle A,B\rangle\rangle$
and showed 
$Gr\varPhi_\mathrm{ Hod}^\sharp(\Phi_{GRT})=\Phi_\mathrm{ KZ}\mod \pi^2$
and
$Gr\varPhi^{(l)\sharp}_{\bold Q_l}(\Phi_{GRT})$ 
$=Gr\Phi^{(l)}_\mathrm{ Ih}$ 
(a series associated with the $l$-adic Ihara associator 
$\Phi^l_\sigma$ in loc.cit).
By \eqref{melon} we get their Berthelot-Ogus counterpart
$$
\varPhi_\mathrm{ crys}^{(p)\sharp}(\Phi_{GRT})=\Phi^p_\mathrm{ KZ}.
$$
\item
In loc.cit \S 6.6.2, we consider the meta-abelian quotient 
$B_{GRT}\in
{\mathcal O}(\underline{GRT}_1)[[A,B]]$ of $\Phi_{GRT}$ 
and showed 
$Gr\varPhi_\mathrm{ Hod}^\sharp(B_{GRT})=
\frac{\Gamma(1-A)\Gamma(1-B)}{\Gamma(1-A-B)}\mod \pi^2$
and
$Gr\varPhi^{(l)\sharp}_{\bold Q_l}(B_{GRT})$ $=GrB^{(l)}(A,B)$
where
$\Gamma(1-z)=
\exp\{\gamma z+\sum\limits_{n=2}^\infty \zeta(n)\frac{z^n}{n}\}$
($\gamma$: Euler constant)
is the classical gamma function
and
$GrB^{(l)}(A,B)$ is
the associated graded quotient of the Ihara's universal power series 
for Jacobi sums \cite{I90}.
By \eqref{melon} we get their Berthelot-Ogus counterpart
$$
\varPhi_\mathrm{ crys}^{(p)\sharp}(B_{GRT})=\prod_{k=1}^\infty
\frac{\Gamma_p(1+p^kA)\Gamma_p(1+p^kB)}{\Gamma_p(1+p^k(A+B))}
$$
where $\Gamma_p(1-z)=
\exp\{\gamma_p z+\sum\limits_{n=2}^\infty
L_p(n,\omega^{1-n})\frac{z^n}{n}\}$
($\gamma_p$: a $p$-adic analogue of Euler constant)
is the Morita's $p$-adic gamma function (\cite{M} Theorem 1).
\end{enumerate}
\end{rem}

Next we discuss a special element of Drinfel'd's 
$\frak{grt}_1=\underset{n\in\bold N}{\oplus}\frak{grt}^n_1$
\footnote{
Actually it is isomorphic to the graded completion of
Ihara's \cite{I90}
stable derivation algebra $\frak D_\centerdot$.
} 
which is the graded Lie algebra of $\underline{GRT}_1$ \cite{Dr} \S 5.
The original Drinfel'd associator $\Phi_\mathrm{ KZ}(A,B)$
satisfies the defining equations of $\underline{M}$ \cite{Dr}
which is similar to but different from 
those of $\underline{GRT}_1$ \cite{Dr}
however we say the following.

\begin{lem}\label{China}
The \lq$-$' part $\Phi^-_\mathrm{ KZ}(A,B)$ 
of the Drinfel'd associator 
satisfies the same equations to those of Proposition \ref{Japan}.
\end{lem}

\begin{pf}
Drinfel'd showed 
$(1,\Phi_\mathrm{ KZ}(\frac{A}{2\pi i},\frac{B}{2\pi i}))
\in\underline{M}_1(\bold C)$ in \cite{Dr}.
By the complex conjugate action we get
$(1,\Phi_\mathrm{ KZ}(\frac{-A}{2\pi i},\frac{-B}{2\pi i}))
\in\underline{M}_1(\bold C)$.
By the torsor structure of $\underline{M}_1(\bold C)$,
there is a unique element $(1,\varphi)\in\underline{GRT}_1(\bold C)$
such that 
$(1,\varphi)\circ(1,\Phi_\mathrm{ KZ}(\frac{-A}{2\pi i},\frac{-B}{2\pi i}))
=(1,\Phi_\mathrm{ KZ}(\frac{A}{2\pi i},\frac{B}{2\pi i}))$,
equivalently 
$$
\Phi_\mathrm{ KZ}\Bigl(\frac{A}{2\pi i},\frac{B}{2\pi i}\Bigr)=
\varphi(A,B)\cdot\Phi_\mathrm{ KZ}\Bigl(\frac{-A}{2\pi i},
\varphi(A,B)^{-1}(\frac{-B}{2\pi i})\varphi(A,B)\Bigr).
$$
By comparing the equations in Lemma \ref{Kazakhstan},
we get $\varphi=\Phi^-_\mathrm{ KZ}(\frac{A}{2\pi i},\frac{B}{2\pi i})$. 
So
$(1,\varphi)=(1,\Phi^-_\mathrm{ KZ}(\frac{A}{2\pi i},\frac{B}{2\pi i}))
\in\underline{GRT}_1$,
from which we get the claim.
\qed
\end{pf}

\begin{rem}\label{Bombay}
By Lemma \ref{China}, 
$(1,\Phi^-_\mathrm{ KZ})\in\underline{GRT}_1(\bold R)$.
By the logarithmic morphism 
$Log:\underline{GRT}_1(\bold R)\to
\frak{grt}_1(\bold R)$, 
Drinfel'd obtained in \cite{Dr} Proposition 6.3
a canonical element $\frak f(\infty)$
in $\frak{grt}_1(\bold R)$,
the image of $(1,\Phi^-_\mathrm{ KZ})$,
with the following presentations:
\begin{equation}\label{Elizabeth}
\frak f(\infty)=\sum\limits_{m\geqslant 3:\text{ odd}}
\frak f(\infty)_m   \  
\text{where}        \  
\frak f(\infty)_m=2\zeta(m)(ad A)^{m-1}(B)+\cdots\in
\frak{grt}^m_1(\bold R).
\end{equation}
On the other hand,
$(1,\Phi^p_\mathrm{ KZ})\in\underline{GRT}_1(\bold Q_p)$
by Proposition \ref{Korea}.
By the logarithmic morphism 
$Log:\underline{GRT}_1(\bold Q_p)\to
\frak{grt}_1(\bold Q_p)$,
we also get a canonical element $\frak f_p$
in $\frak{grt}_1(\bold Q_p)$,
the image of $(1,\Phi^p_\mathrm{ KZ})$,
with the following presentation:
\[
\frak f(p)=\sum\limits_{m\geqslant 3}\frak f(p)_m\qquad
\text{where}\qquad
\frak f(p)_m=\zeta_p(m)(ad A)^{m-1}(B)+\cdots\in
\frak{grt}^m_1(\bold Q_p).
\]
(N.B. $\frak f(\infty)_m=0$ for $m=1$ or even.)
We stress that the summation in $\frak f(\infty)$ 
is taken for odd number
$m$ greater than or equal to 3
while the summation in $\frak f(p)$ is taken for natural number
$m$ greater than or equal to 3.
\end{rem}

\subsection{Racinet's context}\label{Racinet}
We will recall 
Racinet's \cite{R} pro-algebraic torsor $\underline{DMR}_1$
which is defined by two shuffle relations 
and discuss a story for $\underline{DMR}_1$ analogous
to previous section.

MZV $\zeta(k_1,\cdots,k_m)$ ($m,k_1,\cdots,k_m\in \bold N$, $k_m>1$)
satisfies two types of so called shuffle product formula,
expressing a product of two MZV's as a linear combination
of other such values. 
The first type, known as series shuffle product formula 
and for which the easiest example is the relation 
$$
\zeta(k_1)\cdot \zeta(k_2)= \zeta(k_1,k_2) + \zeta(k_2,k_1) + \zeta(k_1+k_2)\;,
$$
is easily obtained from the expression \eqref{gold}.
The second type of shuffle product formula, 
known as integral shuffle product formula, 
come from their iterated integral expressions.
The easiest example is the formula
$$
\zeta(k_1)\cdot\zeta(k_2)=\sum_{i=0}^{k_1-1}\binom{k_2-1+i}{i}\zeta(k_1-i,k_2+i)+\sum_{j=0}^{k_2-1}\binom{k_1-1+j}{j}\zeta(k_2-j,k_1+j).
$$
The double shuffle relations are linear relations combining 
series shuffle product formula and integral shuffle product formula.
The above two formulae give the simplest example 
\begin{align*}
\zeta(k_1,k_2)& + \zeta(k_2,k_1) + \zeta(k_1+k_2) \\
&=\sum_{i=0}^{k_1-1}\binom{k_2-1+i}{i}\zeta(k_1-i,k_2+i)+\sum_{j=0}^{k_2-1}\binom{k_1-1+j}{j}\zeta(k_2-j,k_1+j).
\end{align*}
We have two extended notions of MZV's to non-admissible indices,
the case when the last component $k_m$ is equal to $1$,
which we call series regularization and integral regularization here.
The usual and series regularized MZV's satisfy 
series shuffle product formulae and
the usual and integral regularized MZV's satisfy 
integral shuffle product formulae.
These two kinds of regularized MZV's are related by so called
regularized relations.

Racinet constructed a pro-algebraic variety 
\footnote{
For our convenience slightly we change the definition.}
$\underline{DMR}$
whose set $\underline{DMR}(k)$
of $k$-valued points ($k$: a field of characteristic $0$)
consists of a series  
$\varphi(A,B)\in k\langle\langle A,B\rangle\rangle$
such that $\varphi(A,-B)$ 
satisfies (3.2.1.1) and (3.2.1.2) in \cite{R},
which are relations reflecting the above relations.
Let $\lambda\in k$.
He denoted $\underline{DMR}_\lambda(k)$
to be the subset consisting of the series in $\underline{DMR}(k)$
whose coefficient of $AB$ is equal to $\frac {\lambda^2}{24}$.
He showed in \cite{R} Theorem 1 that
$\underline{DMR}_1(k)$ (resp. $\underline{DMR}_0(k)$)
has a structure of the pro-algebraic torsor 
(resp. pro-algebraic group)
where $\underline{DMR}_0(k)$ acts from the left
(resp. whose group law is defined by the same way to 
the Grothendieck-Teichm\"{u}ller group
$\underline{GRT}_1$).

We embed $\pi_1^{\mathrm{DR}}(X_{\bold Q}:\overrightarrow{01})(k)$
into $k\langle\langle A,B\rangle\rangle$ via $i$ \eqref{navy}
and $\pi_1^\mathrm{ Be}(X(\bold C):\overrightarrow{01})(k)$
into $k\langle\langle A,B\rangle\rangle$ via $j$ \eqref{tan}.
The group automorphism of 
$\pi_1^{\mathrm{DR}}(X_{\bold Q}:\overrightarrow{01})(k)$
sending $A\mapsto A$ and $B\mapsto\varphi^{-1}B\varphi$ induces
an embedding of both
$\underline{DMR}_0$ and 
$\underline{GRT}_1$ into
$\underline{Aut}\pi_1^{\mathrm{DR}}(X_{\bold Q}:\overrightarrow{01})$
as pro-algebraic groups.
And the group isomorphism from 
$\pi_1^{\mathrm{DR}}(X_{\bold Q}:\overrightarrow{01})(k)$ to
$\pi_1^\mathrm{ Be}(X(\bold C):\overrightarrow{01})(k)$
sending $A\mapsto A$ and $B\mapsto\varphi^{-1}B\varphi$ induces
an embedding of both
$\underline{DMR}_1$ and 
$\underline{M}_1$ into
$\underline{Isom}(\pi_1^{\mathrm{DR}}(X_{\bold Q}:\overrightarrow{01}),
\pi_1^\mathrm{ Be}(X(\bold C):\overrightarrow{01}))$
as pro-algebraic torsors.
From now on we regard $\underline{DMR}_0$ and $\underline{GRT}_1$ 
to be sub-pro-algebraic groups of
$\underline{Aut}\pi_1^{\mathrm{DR}}(X_{\bold Q}:\overrightarrow{01})$
and regard $\underline{DMR}_1$ and $\underline{M}_1$ 
to be sub-pro-algebraic torsors of
$\underline{Isom}(\pi_1^{\mathrm{DR}}(X_{\bold Q}:\overrightarrow{01}),
\pi_1^\mathrm{ Be}(X(\bold C):\overrightarrow{01}))$.

We note that 
$\Phi_\mathrm{ KZ}\Bigl(\frac{A}{2\pi i},\frac{B}{2\pi i}\Bigr)$
belongs to $\underline{DMR}_1(\bold C)$,
which gives a morphism 
\begin{equation}\label{nut}
\varPsi_\mathrm{ Hod}:
Spec \frak Z\hookrightarrow \underline{DMR}_1
\end{equation}
analogous to \eqref{peach} and conjectured to be isomorphic \cite{R}.
The paper \cite{BF} is a trial to give a morphism 
\begin{equation}\label{cherry}
\varPsi^{(p)}_\mathrm{ crys}:
Spec \frak Z^{(p)}\hookrightarrow \underline{DMR}_0,
\end{equation}
which is  an analogous map to \eqref{melon}.
We showed double shuffle relations, i.e.
series and integral shuffle product formulae, 
for usual $p$-adic MZV's but not for extended values.
In \cite{FJ} we will give a complete proof to give an embedding 
$\varPsi^{(p)}_\mathrm{ crys}$, namely show 
series and integral shuffle product formulae
and regularized relations for regularized $p$-adic MZV's.

\subsection{Deligne-Goncharov's context}\label{Goncharov}
We will discuss the pro-algebraic bi-torsor of the motivic Galois group
in \cite{DG} 
by recalling the category of mixed Tate motives over $\bold Z$.
This torsor is related with  
the Drinfel'd's torsor (\S \ref{Drinfeld}) in Proposition \ref{Taipei} and
the Racinet's torsor (\S \ref{Racinet}) in Proposition \ref{Taiwan}.
We also explain how the motivic formalism fits into our story in Note \ref{Deligne}:
In Hodge side we give a motivic interpretation of Zagier's dimension conjecture
on MZV's.
We recall how the upper-bounding part of this conjecture follows as a consequence.
In Artin side we see how Ihara's conjecture \cite{I2} on Galois image
is motivically related with Zagier's conjecture via Problem \ref{Belyi}.
In Berthelot-Ogus side
we deduce another motivic way, different from \cite{FJ},
of proving double (series and integral) shuffle relations 
and regularized relations for (extended) $p$-adic multiple zeta values
by using Yamashita's \cite{Y} fiber functor of 
the crystalline realization. 

Let $k$ be a field with characteristic $0$.
Levine \cite{L2} and Voevodsky \cite{Voe} constructed a triangulated 
category of mixed motives over $k$.
Levine \cite{L2} showed an equivalence of these two categories.
This category denoted by $DM(k)_{\bold Q}$ has Tate objects 
$\bold Q(n)$ ($n\in\bold Z$).
Let $DMT(k)_{\bold Q}$ be the triangulated sub-category of 
$DM(k)_{\bold Q}$ generated by $\bold Q(n)$ ($n\in\bold Z$).
Levine \cite{L1} extracted a neutral tannakian $\bold Q$-category 
$MT(k)_{\bold Q}$ of mixed Tate motives over $k$ from $DMT(k)_{\bold Q}$
by taking a heart with respect to a $t$-structure
under the Beilinson-Soul\'{e} vanishing conjecture
which says $gr^\gamma_i K_n(k)=0$ for $n>2i$.
Here LHS is the graded quotient of the algebraic $K$-theory for $k$
with respect to $\gamma$-filtration.

From now on we assume that $k$ is a number field.
In this case the Beilinson-Soul\'{e} vanishing conjecture holds
and we have $MT(k)_{\bold Q}$.
This category satisfies the following expected properties:
Each object $M$ has an increasing filtration of subobjects called 
weight filtration,
$
W:\cdots\subseteq W_{m-1}M\subseteq W_mM\subseteq W_{m+1}M\subseteq\cdots,
$
whose intersection is $0$ and union is $M$.
The quotient $Gr^W_{2m+1}M:=W_{2m+1}M/W_{2m}M$ is trivial and
$Gr^W_{2m}M:=W_{2m}M/W_{2m+1}M$ is a direct sum of 
finite copies of $\bold Q(m)$ for each $m\in\bold Z$.
Morphisms of $MT(k)_{\bold Q}$  are strictly compatible with weight filtration.
The extension group is related to $K$-theory as follows
\begin{equation*}
Ext^i_{MT(k)_{\bold Q}}(\bold Q(0),\bold Q(m))=
\begin{cases}
K_{2m-i}(k)_{\bold Q}  &\text{for  }  i=1,\\
0                      &\text{for  }   i>1.
\end{cases}
\end{equation*}
There are realization fiber functors (\cite{L2} and \cite{H})
corresponding to usual cohomology theories.

Let $S$ be a finite set of finite places of $k$.
Let $\mathcal O_S$ be the ring of $S$-integers in $k$.
Deligne and Goncharov \cite{DG} defined the full subcategory $MT(\mathcal O_S)$
of mixed Tate motives over $\mathcal O_S$,
whose objects are mixed Tate motives $M$ in $MT(k)_{\bold Q}$
such that for each subquotient $E$ of $M$ 
which is an extension of $\bold Q(n)$ by $\bold Q(n+1)$ for $n\in\bold Z$,
the extension class of $E$ in 
$
Ext^1_{MT(k)_{\bold Q}}(\bold Q(n),\bold Q(n+1))=
Ext^1_{MT(k)_{\bold Q}}(\bold Q(0),\bold Q(1))=k^\times_{\bold Q}
$
lies in $\mathcal O_S^\times\otimes\bold Q$.
In this category the following   hold:
\begin{align*}
Ext^1_{MT(\mathcal O_S)}(\bold Q(0),\bold Q(m))&=
\begin{cases}
0                                 &\text{for  }  m<1,\\
\mathcal O_S^\times\otimes\bold Q  &\text{for  } m=1,\\
K_{2m-1}(k)_{\bold Q}             &\text{for  }  m>1,
\end{cases}\\
Ext^2_{MT(\mathcal O_S)}(\bold Q(0),\bold Q(m))&=0.
\end{align*}

Let $\omega_\mathrm{can}:MT(\mathcal O_S)\to Vect_{\bold Q}$
($Vect_{\bold Q}$: the category of $\bold Q$-vector spaces)
be the fiber functor  which sends each motive $M$ to
$\oplus_n Hom(\bold Q(n),Gr^W_{-2n}M)$.
Let $G_\mathrm{can}$ be the motivic Galois group
$\underline{Aut}^\otimes(MT(\bold Z):\omega_\mathrm{can})$.
The action of $G_\mathrm{can}$ on $\omega_\mathrm{can}(\bold Q(1))=\bold Q$
defines a surjection $G_\mathrm{can}\to\bold G_m$ 
and its kernel $U_\mathrm{can}$ is the unipotent radical of $G_\mathrm{can}$.
There is a canonical splitting $\tau:\bold G_m\to G_\mathrm{can}$
which gives a negative grading on the Lie algebra $Lie U_\mathrm{can}$
(consult \cite{De} \S 8 for the full story).
The above computations of $Ext$-groups follows

\begin{prop}\label{Delhi}
The graded Lie algebra $Lie U_\mathrm{can}$ is free and
its degree $n$-part of 
$Lie U_\mathrm{can}^\mathrm{ab}=U_\mathrm{can}^\mathrm{ab}$
is isomorphic to the dual of $Ext^1_{MT(\mathcal O_S)}(\bold Q(0),\bold Q(-n))$.
\end{prop}

\begin{pf}
See \cite{De} \S 8 and \cite{DG} \S 2.
\qed
\end{pf}

Let us restrict in the case of $k=\bold Q$, $S=\emptyset$, 
$\mathcal O_S=\bold Z$.
Let $\omega_*:M\mapsto M_*$ ($*=$ Be, DR) 
be the fiber functor
which associates each mixed Tate motive $M\in MT(\bold Z)$ with 
the underlying vector space of its Betti, De Rham
realization respectively.
For $*,*'=\{\text{Be, DR}\}$
we denote the corresponding tannakian fundamental group 
$\underline{Aut}^\otimes(\mathcal{MT}(\bold Z):\omega_*)$ by $G_*$ and
the corresponding tannakian fundamental torsor
$\underline{Isom}^\otimes(\mathcal{MT}(\bold Z):\omega_*,\omega_{*'})$ 
by $G_{**'}$.
Note that the latter is a $(G_*,G_{*'})$-bi-torsor
and $G_{**}=G_*$.
Let $U_*$ be the sub-pro-algebraic group of $G_*$
whose action on $\omega_*(\bold Q(1))=\bold Q$ is trivial
and $U_{**'}$ be the sub-pro-algebraic torsor of $G_{**'}$
which induces a trivial map from $\omega_*(\bold Q(1))=\bold Q$ 
to $\omega_{*'}(\bold Q(1))=\bold Q$.
Then $U_{**'}$ is a $(U_*,U_{*'})$-bi-torsor and $U_{**}=U_*$.
Since $M_\mathrm{DR}=\omega_\mathrm{can}(M)$ we have 
$G_\mathrm{can}= G_\mathrm{DR}$.
By Proposition \ref{Delhi} the Lie algebra $Lie U_\mathrm{DR}$ of the
unipotent part $U_\mathrm{DR}$ of $G_\mathrm{DR}$ should be
a graded free Lie algebra generated by one element 
in each degree $-m$ ($m\geqslant 3$: odd).

\begin{rem}\label{Allahabad}
In \cite{De} \S 8.12, Deligne constructed a free basis of 
the Lie algebra $Lie U_{\mathrm{DR}}(\bold C)$ as follows:
The infinity Frobenius action $\phi_\infty$ on 
$\omega_{\mathrm{Be}}(M)$ for $M\in\mathcal{MT}(\bold Z)$
determines a point in $G_\mathrm{Be}(\bold Q)$
and by the Hodge comparison isomorphism (the period map)
$\omega_{\mathrm{Be}}(M)\otimes\bold C\simeq
\omega_{\mathrm{DR}}(M)\otimes\bold C$ 
it gives a point in $G_\mathrm{DR}(\bold C)$
which we denote by the same symbol $\phi_\infty$.
Put $\psi_\infty=\phi_\infty\circ\tau(-1)$.
Then $\psi_\infty\in U_{\mathrm{DR}}(\bold C)$
because $\phi_\infty\in G_\mathrm{DR}(\bold C)$ goes to $-1\in\bold G_m$.
The element $\frak e_m(\infty)$  ($m\geqslant 3$: odd) generates freely
the Lie algebra  $Lie U_{\mathrm{DR}}(\bold C)$
where $\frak e_m(\infty)$ means a degree $-m$-part of the image 
$\frak e(\infty)=\sum_m \frak e_m(\infty)$
of $\psi_\infty$ by the logarithmic morphism
$Log:U_{\mathrm{DR}}(\bold C)\to Lie U_{\mathrm{DR}}(\bold C)$.
(N.B. $\frak e_m(\infty)=0$ for $m=1$ or even.) 
\end{rem}

In \cite{DG} \S4 they constructed the motivic fundamental group
$\pi_1^{\mathcal M}(X:\overrightarrow{01})$ with 
$X=\bold P^1\backslash\{0,1,\infty\}$,
which is an ind-object of $MT(\bold Z)$.
This is an affine group $MT(\bold Z)$-scheme (cf. \S \ref{de Rham}),
whose de Rham realization and Betti realization agree with 
the de Rham fundamental group 
$\pi_1^{\mathrm{DR}}(X_{\bold Q}:\overrightarrow{01})$ in \S \ref{de Rham}
and the Betti fundamental group
$\pi_1^{\mathrm{Be}}(X(\bold C):\overrightarrow{01})$ in \S\ref{Betti}
respectively. Namely
$\omega_\mathrm{DR}(\pi_1^{\mathcal M}(X:\overrightarrow{01}))
=\pi_1^{\mathrm{DR}}(X_{\bold Q}:\overrightarrow{01})$ and
$\omega_\mathrm{Be}(\pi_1^{\mathcal M}(X:\overrightarrow{01}))
=\pi_1^{\mathrm{Be}}(X(\bold C):\overrightarrow{01})$.
Since all the structure morphism of $\pi_1^{\mathcal M}(X:\overrightarrow{01})$ 
belong to the set of morphisms of $MT(\bold Z)$ we have 
\begin{equation}\label{Varanasi}
\varphi:G_{**'}\to \underline{Isom}\bigl(
\omega_*(\pi_1^{\mathcal M}(X:\overrightarrow{01})), 
\omega_{*'}(\pi_1^{\mathcal M}(X:\overrightarrow{01}))
\bigr) 
\end{equation}
for $*,*'\in\{\text{Be,DR}\}$.
On this map $\varphi$ the following is one of the basic problems.

\begin{prob}\label{Belyi}
Is $\varphi$ injective?
\end{prob}
This might be said a question which asks a validity of
a unipotent variant of the so-called  `Bely\u{\i}'s theorem' 
in \cite{Be} in the pro-finite setting.
Equivalently this asks if the motivic fundamental group 
$\pi_1^{\mathcal M}(X:\overrightarrow{01})$
is a generator of the tannakian category $MT(\bold Z)$.

The $(U_{\mathrm{DR}},U_{\mathrm{Be}})$-bi-torsor $U_{\mathrm{DR,Be}}$ 
is related to 
the $(\underline{GRT}_1,\underline{GT}_1$)-bi-torsor $\underline{M}_1$
in \S \ref{Drinfeld} as follows.

\begin{prop}\label{Taipei}
$\varphi(U_\mathrm{DR})\subset \underline{GRT}_1$,
$\varphi(U_\mathrm{DR,Be})\subset \underline{M}_1$,
$\varphi(U_\mathrm{Be})\subset \underline{GT}_1$.
\end{prop}

\begin{pf}
By Lemma \ref{Kylgyz}
the automorphism which corresponds to $\phi_\infty\tau(-1)$ 
is described by $(1,\Phi^-_{\mathrm{KZ}})$. 
So the induced homomorphism $\Phi:Lie U_\mathrm{DR}\to 
DerLie \pi_1^{\mathrm{DR}}(X_{\bold Q}:\overrightarrow{01})$
sends $\frak e(\infty)$ to $\frak f(\infty)$ (cf. Remark \ref{Bombay}).
Since $\frak f(\infty)_m$ is the image of the free generator
$\frak e(\infty)_m$,  
the image $\Phi(Lie U_\mathrm{DR})$ must lies in 
$\frak{grt}_1(=Lie \underline{GRT}_1)
\subset DerLie \pi_1^{\mathrm{DR}}(X_{\bold Q}:\overrightarrow{01})$.
So $\varphi(U_\mathrm{DR})$ should lie in $\underline{GRT}_1
\subset \underline{Aut}\pi_1^{\mathrm{DR}}(X_{\bold Q}:\overrightarrow{01})$.

The morphism $\tau(2\pi i)^{-1}p$ lies in $U_\mathrm{DR,Be}(\bold C)$.
Its image $\varphi(\tau(2\pi i)^{-1}p)$,
which is a morphism from 
$\pi_1^{\mathrm{DR}}(X_{\bold Q}:\overrightarrow{01})(\bold C)$ to
$\pi_1^{\mathrm{Be}}(X(\bold C):\overrightarrow{01})(\bold C)$,
is described by 
$(1,\Phi_{\mathrm{KZ}}(\frac{A}{2\pi i},\frac{B}{2\pi i}))$
by Lemma \ref{Vancouver}.
By \cite{Dr} this element belongs to $\underline{M}_1(\bold C)$,
which implies $\varphi(U_\mathrm{DR,Be})\subset \underline{M}_1$.
It is because $\varphi$ is a morphism 
from the $U_\mathrm{DR}$-torsor $U_\mathrm{DR,Be}$
to the $\underline{GRT}_1$-torsor  $\underline{M}_1$,
which sends $\tau(2\pi i)^{-1}p\in U_\mathrm{DR,Be}(\bold C)$ to
$(1,\Phi_{\mathrm{KZ}}(\frac{A}{2\pi i},\frac{B}{2\pi i}))\in
\underline{M}_1(\bold C)$.
\qed
\end{pf}



The $(U_{\mathrm{DR}},U_{\mathrm{Be}})$-bi-torsor $U_{\mathrm{DR,Be}}$ 
is also related to the $\underline{DMR}_0$-torsor $\underline{DMR}_1$
in \S \ref{Racinet} as follows.

\begin{prop}\label{Taiwan}
$\varphi(U_\mathrm{DR})\subset \underline{DMR}_0$,
$\varphi(U_\mathrm{DR,Be})\subset \underline{DMR}_1$.
\end{prop}

\begin{pf}
The proof is similar to the previous proposition.
As for the proof for $\varphi(U_\mathrm{DR})\subset \underline{DMR}_0$,
we use Racinet's result \cite{R} 5.3.2,
which is the same argument to \cite{Dr} Proposition 6.3,
that the Lie algebra $\frak{dmr}_0$
of his $\underline{DMR}_0$
contains the Drinfel'd's element $\frak f(\infty)_m$ \eqref{Elizabeth}.
\qed
\end{pf}

At present we do not know the relationship between
Drinfel'd's  $\underline{GRT}_1$ and Racinet's $\underline{DMR}_0$
although it might be expected that they are equal.
However
the recent Terasoma's work (to appear) might suggests a direction
$\underline{DMR}_0\subset\underline{GRT}_1$.

Below we explain how the motivic formalism fits in the scheme of our subject
in the previous section.
\begin{note}\label{Deligne}
\begin{enumerate}
\item
In Hodge side: as is shown in the proof of Proposition \ref{Varanasi}
the pair $(1,\Phi_{\mathrm{KZ}}(\frac{A}{2\pi i},\frac{B}{2\pi i}))$
is the image of $\tau(2\pi i)^{-1}p$ by $\varphi$, that is,
$(1,\Phi_{\mathrm{KZ}}(\frac{A}{2\pi i},\frac{B}{2\pi i}))\in
\varphi(U_\mathrm{DR,Be})(\bold C)$.
Since $(2\pi i)^{-(k_1+\cdots+k_m)}\zeta(k_1,\dots,k_m)$
($m,k_i\in\bold N$,$k_m>1$) appears among each coefficient of 
$\Phi_{\mathrm{KZ}}(\frac{A}{2\pi i},\frac{B}{2\pi i})$,
it gives rise a morphism
\begin{equation}\label{Agra}
\Gamma_\mathrm{Hod}:Spec\mathcal{Z}\hookrightarrow\varphi(U_\mathrm{DR,Be})
\end{equation}
analogous to \eqref{peach} and \eqref{nut}.
There is a conjecture on the dimension of the vector space of MZV's at each weight
which is called Zagier conjecture (see \cite{Za} and also \cite{F0})
and partly proved (for upper-bounding part) by Terasoma \cite{T}.
By using the embedding \eqref{Agra}  and Proposition \ref{Delhi},
Deligne and Goncharov also get a partial (upper-bounding part) proof of
Zagier conjecture in \cite{DG}.
We also get to know that 
to say Zagier conjecture holds is equivalent to say
the surjectivity of $\Gamma_\mathrm{Hod}$ and
the injectivity of $\varphi$ (i.e. the validity of Problem \ref{Belyi}).
\item
In Artin side: we consider the absolute Galois group 
$Gal(\overline{\bold Q}/\bold Q)$ action on the $l$-adic ($l$: a prime)
\'{e}tale fundamental group 
$\varphi_1:Gal(\overline{\bold Q}/\bold Q)\to\underline{Aut}
\pi_1^{l,\text{\'{e}t}}(X_{\overline{\bold Q}}:\overrightarrow{01})$
and the $l$-adic Galois image pro-algebraic group
$\underline{Gal}^{(l)}_{\bold Q_l}$
defined to be the Zariski closure of 
$\varphi_1\bigl(Gal(\overline{\bold Q}/\bold Q(\mu_{l^\infty})\bigr)$ 
in \cite{F0} \S 4.
The category $MT(\bold Z)$ has
a fiber functor $\omega_{l,\text{\'{e}t}}:M\mapsto M_{l,\text{\'{e}t}}$
of the $l$-adic \'{e}tale realization 
which associates each motives $M$ with its underlying $\bold Q_l$-vector space
of its $l$-adic \'{e}tale realization.
The absolute Galois group $Gal(\overline{\bold Q}/\bold Q)$ acts functorially on 
$\omega_{l,\text{\'{e}t}}(M)$ and hence it induces a morphism 
$Gal(\overline{\bold Q}/\bold Q)\to G_{l,\text{\'{e}t}}(\bold Q_l)$.
Hereafter we fix an embedding $\overline{\bold Q}\hookrightarrow\bold C$.
Then there is the Artin comparison isomorphism 
$\omega_{l,\text{\'{e}t}}(M)\simeq \omega_\mathrm{Be}(M)\otimes\bold Q_l$
functorial with respect to $M$,
which gives 
$G_{l,\text{\'{e}t}}=G_\mathrm{Be}\times\bold Q_l$.
We might have a morphism $\varphi_2:Gal(\overline{\bold Q}/\bold Q)\to 
G_\mathrm{Be}(\bold Q_l)$.
Because of
$\pi_1^{l,\text{\'{e}t}}(X_{\overline{\bold Q}}:\overrightarrow{01})=
\omega_{l,\text{\'{e}t}}(\pi_1^{\mathcal M}(X:\overrightarrow{01}))=
\omega_\mathrm{Be}(\pi_1^{\mathcal M}(X:\overrightarrow{01}))\times\bold Q_l$,
the map $\varphi_1$ might factor through $G_\mathrm{Be}$-action on 
$\omega_{Be}(\pi_1^{\mathcal M}(X:\overrightarrow{01}))$, that means
$\varphi_1=\varphi\circ\varphi_2$
(consult \cite{De} \S 8 for full story).
In loc.cit. \S 8.14 it is also shown that the image of $\varphi_2$ is open 
in the topology of $G_\mathrm{Be}(\bold Q_l)$.
So by restricting into unipotent parts we get an isomorphism
\begin{equation}
\Gamma_{l,\text{\'{e}t}}:\underline{Gal}^{(l)}_{\bold Q_l}
\overset{=}{\to}
\varphi(U_\mathrm{Be})\times\bold Q_l
\end{equation}
which might be said an analogue of \eqref{apple}.
This is the way how the common $\bold Q$-structure of the $l$-adic Galois image 
pro-algebraic group for all prime $l$ is detected in \cite{DG} Remark 6.13.
The graded Lie algebra $\frak g^{(l)}_\centerdot$ (cf.\cite{I2} and also \cite{F0})
associated with $\underline{Gal}^{(l)}_{\bold Q_l}$ 
is conjectured by Ihara (loc.cit.) to be a free Lie algebra generated by 
one element in each degree 3,5,7,9,..
As in Hodge case we can say that this Ihara's conjecture is equivalent to
the injectivity of $\varphi$ (i.e. the validity of Problem \ref{Belyi}).
\item
In Berthelot-Ogus side: Yamashita will construct in \cite{Y} the fiber functor 
$\omega_{p,\mathrm{crys}}:M\to M_{p,\mathrm{crys}}$
of the crystalline realization of $MT(\bold Z)$
which is compatible with $\omega_\mathrm{DR}$.
It associates each mixed Tate motive $M$ with its
underlying $\bold Q_l$-vector space of its crystalline realization
by using Fontaine's functor
to give a conjectured dimension bounding for $Z^{(p)}_w$.
It admits
a functorial crystalline Frobenius action $\phi_{p,\mathrm{crys}}$
on $\omega_{p,\mathrm{crys}}(M)$ and a functorial comparison isomorphism 
$\omega_{p,\mathrm{crys}}(M)=\omega_\mathrm{DR}(M)\otimes \bold Q_p$
for $MT(\bold Z)$,
which gives a point in $G_\mathrm{DR}(\bold Q_p)$
denoted by the same symbol $\phi_{p,\mathrm{crys}}$.
Since $\phi_{p,\mathrm{crys}}$ goes to $\frac{1}{p}\in\bold G_m$,
the element $\psi_p=\phi_{p,\mathrm{crys}}\circ\tau(p)$
must lie in $U_{\mathrm{DR}}(\bold Q_p)$.
As is similar to Remark \ref{Bombay} and Remark \ref{Allahabad}
we also get in the $p$-adic setting elements $e(p)_m$ ($m\geqslant 3$)
which is on a degree $-m$-part of $Lie U_\mathrm{DR}(\bold Q_p)$.
The action of $\psi_p$ on $\pi_1^{\mathrm{DR}}(X_{\bold Q}:\overrightarrow{01})$
is described as $A\mapsto A$ and 
$B\mapsto \Phi^{p}_{\mathrm{De}}(A,B)^{-1}B\Phi^p_{\mathrm{De}}(A,B)$
by Lemma \ref{Ireland}.
Therefore 
$(1,\Phi^{p}_{\mathrm{De}})\in\varphi(U_\mathrm{DR}(\bold Q_p))$.
By Proposition \ref{Taiwan}, 
$\Phi^p_{\mathrm{De}}\in\underline{DMR}_0(\bold Q_p)$.
Since Deligne's $p$-adic MZV's $\zeta^{\mathrm{De}}_p(k_1,\cdots,k_m)$
are coefficients of $\Phi^p_{\mathrm{De}}(A,B)$,
they must satisfy 
double shuffle relations. 
By the same argument to Proposition \ref{Japan},
we get $(1,\Phi^p_{\mathrm{KZ}})\in\varphi(U_\mathrm{DR}(\bold Q_p))$.
It follows that we also get an embedding
\begin{equation}\label{desk}
\Gamma^{(p)}_\mathrm{ crys}:
Spec \frak Z^{(p)}\hookrightarrow\varphi(U_\mathrm{DR}).
\end{equation}
Combining \eqref{desk} with Proposition \ref{Taiwan}, 
we get the embedding \eqref{cherry}.
Therefore our $p$-adic MZV's $\zeta_p(k_1,\cdots,k_m)$
must satisfy (series and integral) double shuffle relations and 
regularization relations, which are defining relations of $\underline{DMR}_0$.
This is a motivic proof of double shuffle relations for
$p$-adic MZV's,
which was a project posed in \cite{De2}.
\end{enumerate}
\end{note}


\end{document}